\documentclass{amsart}
\usepackage[english]{babel}

\usepackage{amstext}
\usepackage{amssymb}
\usepackage{amsmath}
\usepackage{amsrefs}
\usepackage{mathrsfs}
\usepackage{graphicx}
\usepackage{array}
\usepackage{mathrsfs}
\usepackage{verbatim}
\usepackage{braket}
\usepackage[margin=1.2in]{geometry}
\usepackage{colortbl}
\usepackage{hyperref}

\numberwithin{equation}{section}

\theoremstyle{plain}

\newtheorem{example}{Example}[section]
\newtheorem*{theorem*}{Theorem}

\newtheorem*{prop*}{Proposition}

\theoremstyle{definition}

\theoremstyle{remark}
\newtheorem*{remark}{Remark}

\DeclareMathOperator{\tr}{tr}

\newcommand{\ud}{\,\mathrm{d}}
\newcommand{\RR}{\mathbb{R}}

\newcommand{\NN}{\mathbb{N}}

\newcommand{\bd}[1]{\boldsymbol{#1}}
\newcommand{\I}{\mathrm{i}}
\newcommand{\abs}[1]{\lvert#1\rvert}
\newcommand{\mc}[1]{\mathcal{#1}}
\newcommand{\average}[1]{\langle#1\rangle}
\renewcommand{\Im}{\mathfrak{Im}}

\usepackage{color}
\newcommand{\wh}[1]{\widehat{#1}}
\newcommand{\wt}[1]{\widetilde{#1}}
\newcommand{\wb}[1]{\widebar{#1}}
\newcommand{\norm}[1]{\lVert#1\rVert}

\newcommand{\FGA}{\mathrm{FGA}}

\IfFileExists{mathabx.sty}%
  {\DeclareFontFamily{U}{mathx}{\hyphenchar\font45}%
   \DeclareFontShape{U}{mathx}{m}{n}{<->mathx10}{}%
   \DeclareSymbolFont{mathx}{U}{mathx}{m}{n}%
   \DeclareFontSubstitution{U}{mathx}{m}{n}%
   \DeclareMathAccent{\widebar}{0}{mathx}{"73}%
}{%
  \PackageWarning{mathabx}{%
    Package mathabx not available, therefore\MessageBreak substituting
    widebar with overline\MessageBreak }%
  \newcommand{\widebar}[1]{\wb{#1}}%
}

\newcommand{\veps}{\varepsilon}

\begin{document}

\title[Gauge-invariant frozen Gaussian approximation
method]{Gauge-invariant frozen Gaussian approximation method for the
  Schr\"odinger equation with periodic potentials}

\date{\today}

\author{Ricardo Delgadillo}
\address{Department of Mathematics, University of California, Santa Barbara,
CA 93106}
\email{rdelgadi@math.ucsb.edu}

\author{Jianfeng Lu}
\address{Departments of Mathematics, Physics, and
  Chemisty, Duke University, Box 90320, Durham, NC 27708}
\email{jianfeng@math.duke.edu}

\author{Xu Yang}
\address{Department of Mathematics, University of California, Santa Barbara,
CA 93106}
\email{xuyang@math.ucsb.edu}

\thanks{R.D. and X.Y. were partially supported by the NSF grants
  DMS-1418936 and DMS-1107291: NSF Research Network in Mathematical
  Sciences ``Kinetic description of emerging challenges in multiscale
  problems of natural science''.  They also acknowledge support from
  the Center for Scientific Computing at the CNSI and MRL: an NSF
  MRSEC (DMR-1121053) and NSF CNS-0960316.  The work of J.L.~was
  supported in part by the Alfred P.~Sloan Foundation and the National
  Science Foundation under award DMS-1312659 and DMS-1454939. X.Y. was
  also partially supported by the Regents Junior Faculty Fellowship
  and Hellman Family Foundation Faculty Fellowship of University of
  California, Santa Barbara. }

\maketitle

\begin{abstract}
  We develop a gauge-invariant frozen Gaussian approximation (GIFGA)
  method for the linear Schr\"odinger equation (LSE) with periodic
  potentials in the semiclassical regime. The method generalizes the
  Herman-Kluk propagator for LSE to the case with periodic media.  It
  provides an efficient computational tool based on asymptotic
  analysis on phase space and Bloch waves to capture the
  high-frequency oscillations of the solution. Compared to geometric
  optics and Gaussian beam methods, GIFGA works in both scenarios of
  caustics and beam spreading. Moreover, it is invariant with respect
  to the gauge choice of the Bloch eigenfunctions, and thus avoids the
  numerical difficulty of computing gauge-dependent Berry phase. We
  numerically test the method by several one-dimensional examples, in
  particular, the first order convergence is validated, which agrees
  with our companion analysis paper [Delgadillo, Lu and Yang,
  arXiv:1504.08051].
\end{abstract}

\section{Introduction}

The focus of this work is to develop efficient numerical methods for
solving the following semiclassical Schr\"{o}dinger equation whose
potential term consists of a (highly oscillatory) microscopic periodic
potential and a macroscopic smooth potential,
\begin{equation}\label{eq:schrodinger}
\I\veps\dfrac{\partial\psi^{\veps}}{\partial t}=-\dfrac{\veps^{2}}{2}\Delta\psi^{\veps}+V_{\Gamma}\left(\dfrac{\bd{x}}{\veps}\right)
\psi^{\veps}+U(\bd{x})\psi^{\veps},\hspace{1cm}\bd{x}\in\mathbb{R}^{d}.
\end{equation}
Here $\psi^{\veps}(t,\bd{x})$ is the wave function and $\veps\ll 1$ is
an effective Planck constant.  The equation \eqref{eq:schrodinger}
can be viewed as a model for electron dynamics in crystal under the
one-particle approximation. The periodic lattice potential $V_\Gamma$
is generated by the ionic cores and electrons in the crystal, and
hence periodic with respect to the lattice $\mathbb{L}$ with unit cell
$\Gamma:=[-\pi,\pi)^{d}$. In \eqref{eq:schrodinger}, $U$ is a smooth external macroscopic
potential, which counts for e.g., external electric field.

Direct numerical simulation of \eqref{eq:schrodinger} is prohibitively
expensive due to the small parameter $\veps$ in the semiclassical
regime. In order to accurately capture the small scale features caused
by $V_{\Gamma}$, a mesh size of order $o(\veps)$ is usually required
in time and space, e.g., in the standard time-splitting spectral
method \cite{BaJiMa:02}. If only physical observables ({\it e.g.},
density, flux and energy) are needed, one can relax the time step
requirement to $\mathcal{O}(1)$ with a coarser mesh size of
$\mathcal{O}(\veps)$ using the Bloch decomposition based
time-splitting spectral method as proposed in
\cite{HuJiMaSp:07,HuJiMaSp:08,HuJiMaSp:09}. However, computation of
the solution $\psi^{\veps}$ to \eqref{eq:schrodinger} is still very
expensive for $\veps\ll1$, especially in high dimensions. For this
reason, alternative approaches based on asymptotic analysis have been
developed, among which, the geometric optics (GO) approach is based on the
WKB ansatz under the adiabatic approximation,
\[\psi^{\veps}(t,\bd{x})=a(t,\bd{x})u_n\left(\nabla_{\bd{x}}S,\frac{\bd{x}}{\veps}\right)e^{\I
  S(t,\bd{x})/\veps}.\]
Here $u_n(\bd{\xi},\bd{x})$ is the Bloch eigenfunction
normalized for each $\bd{\xi}\in\Gamma^{\ast}:=[0,1)^{d}$:
\begin{equation}\label{eq:normal}
  \int_{\Gamma}|u_{n}(\bd{\xi},\bd{x})|^{2} \ud \bd{x}=1,
\end{equation}
which corresponds to the $n$-th energy band $E_n(\bd{\xi})$ (see e.g.,
\cite{BeLiPa:78}):
\begin{equation}\label{eq:eigen}
H_{\bd{\xi}} u_n(\bd{\xi},\bd{x})=E_n(\bd{\xi}) u_n(\bd{\xi},\bd{x}),
\end{equation}
with the Bloch Hamiltonian
\begin{equation}
H_{\bd{\xi}} :=\dfrac{1}{2}(-\I\partial_{\bd{x}}+\bd{\xi})^{2}+V_{\Gamma}(\bd{x}), \label{eq:Hamiltonian}
\end{equation}
and periodic boundary conditions on $\Gamma$.

Then GO solves $S(t,\bd{x})$ as the solution
to an eikonal equation and $\rho(t,\bd{x})=|a(t,\bd{x})|^{2}$ given by a
transport equation:
\begin{align}\label{eq:Eikonal1}
& \partial_{t}S+E_n(\nabla_{\bd{x}}S)+U(\bd{x})=0, \\
\label{eq:Eikonal2}
& \partial_{t}\rho+\nabla_{\bd{x}}\cdot\bigl(\rho\nabla_{\bd{\xi}}E_n(\nabla_{\bd{x}}S)\bigr)=0.
\end{align}
While this method is $\veps$-independent, it breaks down at caustics
where the Hamilton-Jacobi equation \eqref{eq:Eikonal1} develops
singularities.

The Gaussian beam method (GBM) was proposed in \cite{He:81,He:91} to
overcome this drawback at caustics, with some recent developments
\cite{DiGuRa:06,JiWuYa:08,JiWuYaHu:10,JiMaSp:11,YiZh:11,WuHuJiYi:12,JeJi:14},
which in particular extends the method to periodic media. GBM is based on the single beam solution, which has a similar
form as the WKB ansatz
\[ \psi^\veps(t,\bd{x})=a(t,\bd{y})u_n\left(\nabla_{\bd{x}}S,\frac{\bd{x}}{\veps}\right)e^{\I \widetilde{S}(t,\bd{x},\bd{y})/\veps}.
\]
The difference lies in that GBM uses a {complex} phase function,
\begin{equation}\label{eq:GBphs}
  \widetilde{S}(t,\bd{x},\bd{y})=S(t,\bd{y})+\bd{p}(t,\bd{y})\cdot(\bd{x}-\bd{y})+
  \frac{1}{2}(\bd{x}-\bd{y})\cdot
  M(t,\bd{y})(\bd{x}-\bd{y}),
\end{equation}
where $S\in\RR,\;\bd{p}\in\RR^d,\;M\in\mathbb{C}^{d\times d}$. The
imaginary part of $M$ is chosen to be positive definite so that the
solution decays exponentially away from $\bd{x}=\bd{y}$ as a Gaussian,
where $\bd{y}$ is called the beam center. If the initial wave is not
in a form of single beam, one can approximate it by a superposition of
Gaussian beams. The validity of this construction at caustics was
analyzed in \cite{DiGuRa:06}.

The accuracy of GBM relies on the truncation
error of the Taylor expansion of $\widetilde{S}$ around the beam center
$\bd{y}$ up to the quadratic term, and thus it loses accuracy when the
width of the beam becomes large, i.e., when the imaginary part of
$M(t, \bd{y})$ in \eqref{eq:GBphs} becomes small so that the Gaussian
function is no longer localized. This happens for example when the
solution of the Schr\"odinger equation spreads (the opposite situation
of forming caustics). This is a severe problem in general, as shown in
\cite{LuYa:11, MoRu:10, QiYi2:10}. One can overcome the problem of
spreading of beams by doing reinitialization once in a while, see
\cites{QiYi:10,QiYi2:10}, however, this increases the computational complexity
especially when beams spread quickly.

In the setting of semiclassical Schr\"odinger equations with periodic
potential, another challenge for asymptotics methods, not emphasized enough in the
literature though, comes from the gauge
freedom in \eqref{eq:eigen}. That is, for any Bloch eigenfunction
$u_n(\bd{\xi},\bd{x})$, $u_n(\bd{\xi},\bd{x})e^{\I\phi(\bd{\xi})}$
also solves \eqref{eq:eigen} for any arbitrary phase function
$\phi(\bd{\xi})$.  In particular, when one solves the Bloch waves
numerically from \eqref{eq:eigen} for different $\bd{\xi}$, it is very
difficult, if not impossible, to make sure that the phase depends
smoothly on $\bd{\xi}$. The arbitrariness creates difficulty when
one needs to get the eigenfunctions off numerical grids by
interpolation, e.g., in the Gaussian beam method \cite{DiGuRa:06}.

In this paper, we develop a gauge-invariant frozen Gaussian
approximation (GIFGA) method for the Schr\"odinger equation with periodic
potentials.  This method generalizes the Herman-Kluk propagator
\cite{HeKl:84} by including Bloch waves in the integral
representation. It provides an efficient
computational tool based on asymptotic analysis on phase plane, with a
first order accuracy established in our companion analysis paper
\cite{DeLuYa:analysis}.  It inherits the merits of the frozen Gaussian
approximation studied in \cite{LuYa:11,LuYa:12,LuYa2:12}, which works
in both scenarios of caustics and beam spreading. The formulation is
also invariant with respect to the gauge choice of the Bloch
eigenfunctions.  In particular, we avoid the numerical computation of the
Berry phase, which causes difficulty since it depends on the
derivatives of Bloch eigenfunctions with respect to crystal momentum,
and is hence not always well-defined if an arbitrary gauge choice was
made.
This is achieved by using a trick inspired by the work of Vanderbilt
and King-Smith \cite{PhysRevB.47.1651} in the context of modern theory
of polarization. The details will be explained in
Section~\ref{nogauge}, see in particular,
\eqref{eq:SA}--\eqref{eq:Fn}.

The rest of the paper is organized as follows: In Section~\ref{sec:Bloch}, we will introduce the GIFGA method. In Section~\ref{sec:algorithm}, we briefly describe how to numerically compute Bloch eigenvalues and eigenfunctions. We also describe how to numerically implement the GIFGA method described in Section~\ref{sec:Bloch}. Section~\ref{sec:example} presents numerical evidence supporting the initial decomposition described in Section~\ref{sec:Bloch} along with examples confirming our analytical results in \cite{DeLuYa:analysis}. The last two examples in Section~\ref{sec:example} provides the numerical performance of GIFGA. We make some concluding remarks in Section~\ref{sec:conclusion}.

\section{Formulation of the frozen Gaussian
  approximation}\label{sec:Bloch}
This section is devoted to the development of the gauge-invariant
frozen Gaussian approximation (GIFGA) in periodic media based on Bloch
decomposition. We first recall the Bloch decomposition for
Schr\"odinger operators with a periodic potential. The Bloch waves
will be used to capture the high-frequency oscillatory structure of the
solution given by GIFGA. After stating the asymptotic solution,
the formulation of which is gauge-invariant, we recall some analytical
results on the convergence of GIFGA.


\subsection{The Bloch decomposition}

Recall that the potential $V_{\Gamma}(\bd{x})$ in
\eqref{eq:schrodinger} is smooth and periodic with respect to the
lattice $\mathbb{L}$ with unit cell $\Gamma=[-\pi,\pi)^{d}$. The unit
cell of the reciprocal lattice, known as the first Brillouin zone, is
then given by $\Gamma^{\ast}=[0,1)^{d}$.


The eigenvalues of the self-adjoint Bloch Hamiltonian $H_{\bd{\xi}}$,
defined in \eqref{eq:Hamiltonian} on $L^2(\Gamma)$ are real and ordered
increasingly (counting multiplicity) as
\begin{equation}
E_{1}(\bd{\xi})\leq E_{2}(\bd{\xi})\leq\cdots\leq E_{n}(\bd{\xi})\leq\cdots,\hspace{1cm}n\in\mathbb{N}
\end{equation}
for each $\bd{\xi}\in\Gamma^{\ast}$. Furthermore, the eigenfunctions
$\left\{u_{n}(\bd{\xi},\bd{x})\right\}_{n=1}^{\infty}$ for each
$\bd{\xi}\in\Gamma^{\ast}$, known as the Bloch waves, form an
orthonormal basis of $L^{2}(\Gamma)$ \cite{BeLiPa:78}.

We extend $u_{n}(\bd{\xi},\bd{x})$ periodically with respect to
$\bd{x}$ so that it is defined on all of $\mathbb{R}^{d}$, and then the Bloch decomposition is given by, $\forall f\in L^2(\mathbb{R}^d)$,
\begin{equation}\label{eq:Band Decomposition}  f(\bd{x})=\dfrac{1}{(2\pi)^{d/2}}\displaystyle\sum_{n=1}^{\infty}\displaystyle\int_{\Gamma^{\ast}}u_{n}(\bd{\xi},\bd{x})e^{\I\bd{\xi}\cdot\bd{x}}(\mathcal{B}f)_{n}
  (\bd{\xi})\ud\bd{\xi},
\end{equation}
where the Bloch transform $\mathcal{B}:L^{2}(\mathbb{R}^{d})\rightarrow L^{2}(\Gamma^{\ast})^{\mathbb{N}}$ is given by
\begin{equation}
(\mathcal{B}f)_{n}(\bd{\xi})=\dfrac{1}{(2\pi)^{d/2}}\int_{\mathbb{R}^{d}}\overline{u}_{n}(\bd{\xi},\bd{y})e^{-\I\bd{\xi}\cdot\bd{y}}f(\bd{y})\ud\bd{y}.
\end{equation}
As an analog to the Parseval's identity, it holds
\begin{equation}
\int_{\mathbb{R}^{d}}|f(\bd{x})|^{2}\ud\bd{x}=\sum_{n=1}^{\infty}\int_{\Gamma^{\ast}}
\bigl\lvert (\mathcal{B}f)_{n}(\bd{\xi})\bigr\rvert^{2}\ud\bd{\xi}.
\end{equation}

We denote $\Omega$ the phase space corresponding to one band
\begin{equation}
\Omega :=\mathbb{R}^{d}\times\Gamma^{\ast}=\{(\bd{x},\bd{\xi})\mid\bd{x}\in\mathbb{R}^{d},\bd{\xi}\in\Gamma^{\ast}\}.
\end{equation}

Let us
define the Berry phase, which will be used later, as
\begin{equation}\label{eq:berry}
\mathcal{A}_{n}(\bd{\xi})=\langle u_{n}(\bd{\xi},\cdot)|\I\nabla_{\bd{\xi}}u_{n}(\bd{\xi},\cdot)\rangle_{L^2(\Gamma)}.
\end{equation}
Here we have used the Dirac bra-ket notation $\langle\cdot|\cdot\rangle$ in quantum mechanics, i.e.,
\begin{equation*}
\langle f | g \rangle_{L^2(\Omega)}= \int_\Omega\bar{f}g\ud y, \quad \text{and} \quad
\langle f | g \rangle= \int_{\mathbb{R}^d}\bar{f}g\ud y,
\end{equation*}
where $\bar{f}$ is the complex conjugate of $f$. Note that the eigenvalue equation \eqref{eq:eigen} and its normalization only define $u_{n}(\bd{\xi},\cdot)$ up to a unit complex number, in particular, for any function $\phi$ periodic in $\Gamma^{\ast}$
\begin{equation}
\widetilde{u}_{n}(\bd{\xi},\bd{x})=e^{\I\phi(\bd{\xi})}u_{n}(\bd{\xi},\bd{x})
\end{equation}
also provides a set of Bloch waves. This is known as the gauge freedom
for Bloch waves. It is known that (see {\it e.g.},~\cite{ReedSimon4})
we can choose $\phi$ such that $\widetilde{u}_{n}(\bd{\xi},\bd{x})$ is
smooth in $\bd{\xi}$, and then the definition \eqref{eq:berry} makes
sense. However, different gauge choice might give different values of
$\mathcal{A}_{n}(\bd{\xi})$, and it is also difficult in numerical
diagonalization of the Bloch waves to make sure that the phase
dependence is smooth. We will come back to this delicacy in the
development of numerical algorithms. Note that from the normalization
condition \eqref{eq:normal}, $\mathcal{A}_{n}(\bd{\xi})$ is always a real number.

\subsection{Formulation}

We denote $G_{\bd{q},\bd{p}}^{\veps}$ the semiclassical Gaussian
function localized in the phase space at $(\bd{q},\bd{p})$:
\begin{equation}\label{eq:Gaussian}
G_{\bd{q},\bd{p}}^{\veps}(\bd{x})=\exp\bigl(-|\bd{x}-\bd{q}|^{2}/(2\veps)+\I\bd{p}\cdot(\bd{x}-\bd{q})/\veps\bigr).
\end{equation}

The frozen Gaussian approximation (FGA) solution to \eqref{eq:schrodinger} with the initial condition $\psi_{0}$ is approximated by \cite{DeLuYa:analysis},
\begin{multline}\label{eq:FGA}
  \psi_{\FGA}^{\veps}(t,\bd{x})=\dfrac{1}{(2\pi\veps)^{3d/2}}
  \sum_{n=1}^{\infty} \int_{\Omega}a_{n}(t,\bd{q},\bd{p})u_{n}(\bd{P}_{n},\bd{x}/\veps)G_{\bd{Q}_{n},\bd{P}_{n}}^{\veps}(\bd{x})e^{\I S_{n}(t,\bd{q},\bd{p})/\veps}  \\
 \times \braket{G_{\bd{q},\bd{p}}^{\veps}u_{n}(\bd{p},\cdot/\veps)|\psi_{0}}\ud\bd{q}\ud\bd{p}.
\end{multline}
The right hand side of \eqref{eq:FGA} sums over all the Bloch
bands. For each $n$,
$\left(\bd{Q}_{n}(t,\bd{q},\bd{p}),\bd{P}_n(t,\bd{q},\bd{p})\right)$
solves the equation of motion given by the classical Hamiltonian
$h_{n}(\bd{q},\bd{p})=E_{n}(\bd{p})+U(\bd{q})$:
\begin{equation}\label{eq:Hamiltonian Flow}
\left\{\begin{matrix}
\dfrac{d\bd{Q}_{n}}{dt}=\nabla E_{n}(\bd{P}_{n}),\\[1em]
\dfrac{d\bd{P}_{n}}{dt}=-\nabla U(\bd{Q}_{n}),\\
\end{matrix}\right.
\end{equation}
with the initial conditions $\bd{Q}_{n}(0,\bd{q},\bd{p})=\bd{q}$ and $\bd{P}_{n}(0,\bd{q},\bd{p})=\bd{p}$. For simplicity, we shall omit the subscripts of gradient whenever it does not cause any confusion.

In \eqref{eq:FGA}, $S_{n}(t,\bd{q},\bd{p})$ is the action associated with the Hamiltonian dynamics \eqref{eq:Hamiltonian Flow}, given by the evolution equation
\begin{equation}
\dfrac{dS_{n}}{dt}=\bd{P}_{n}\cdot\nabla_{\bd{P}}h_{n}(\bd{Q}_{n},\bd{P}_{n})-h_{n}(\bd{Q}_{n},\bd{P}_{n}),
\end{equation}
with the initial condition $S_{n}(0,\bd{q},\bd{p})=0$. The function $a_{n}(t,\bd{q},\bd{p})$ gives the amplitude of the Gaussian function at time $t$. With  the short hand notations
\begin{equation}
\partial_{\bd{z}}=\partial_{\bd{q}}-\I\partial_{\bd{p}}\hspace{1cm}Z_{n}=\partial_{\bd{z}}\left(\bd{Q}_{n}+\I\bd{P}_{n}\right),
\end{equation}
the evolution equation for $a_{n}$ is given by
\begin{equation}\label{eq:Aevolution}
\dfrac{da_{n}}{dt}=-\I a_{n}\mathcal{A}_{n}(\bd{P}_{n})\cdot\nabla U(\bd{Q}_{n})+\dfrac{1}{2}a_{n}\text{tr}\left(\partial_{z}\bd{P}_{n}\nabla^{2}E_{n}(\bd{P}_{n})Z_{n}^{-1}\right)-\dfrac{\I}{2}a_{n}\text{tr}\left(\partial_{z}\bd{Q}_{n}\nabla^{2}U_{n}(\bd{Q}_{n})Z_{n}^{-1}\right)
\end{equation}
with initial condition $a_{n}(0,\bd{q},\bd{p})=2^{d/2}$ for each
$(\bd{q},\bd{p})$. Recall that $\mathcal{A}_{n}(\bd{\xi})$ is the
Berry phase of the $n$-th Bloch band given in \eqref{eq:berry}.

\subsection{Gauge-Invariant Integrator}\label{nogauge}
The gauge freedom of the eigenfunction $u_{n}(\bd{\xi},\bd{x})$ of
\eqref{eq:Hamiltonian} causes problems for numerical computation. In
particular, different choice of gauge may lead to different numerical
results for the Berry phase term
$\mathcal{A}_{n}(\bd{\xi})=\langle
u_{n}(\bd{\xi},\bd{x})|\I\nabla_{\bd{\xi}}u_{n}(\bd{\xi},\bd{x})\rangle$,
and hence different $\psi_{\FGA}^{\veps}$ which is artificial. It is
desirable hence to design an algorithm that is manifestly independent
of the gauge. The key is to avoid direct computation of the the Berry
phase and so to avoid the the computation of the momentum-gradient of
$u_n$.

First, we separate the dependence of $a_n$ on $\mc{A}_n$ in the evolution
equation \eqref{eq:Aevolution}. For this, we define $S_{n}^{\mc{A}}$ the phase contribution due to the Berry phase
term
\begin{equation}\label{eq:defSA}
  S_n^{\mc{A}}(t, \bd{q}, \bd{p}) = \int_0^t \mc{A}_n(\bd{P}_n)\cdot
  \nabla U(\bd{Q}_n) \ud s.
\end{equation}
Let
\begin{equation*}
  b_n(t, \bd{q}, \bd{p}) = a_n(t, \bd{q}, \bd{p}) \exp(
  \I S_n^{\mc{A}}(t, \bd{q}, \bd{p})),
\end{equation*}
then it solves
\begin{equation}\label{eq:bamplitude}
  \frac{\ud b_n}{\ud t}= \frac{1}{2}b_n \tr\Bigl(
  \partial_{\bd{z}}\bd{P}_{n}
  \nabla^2 E_n(\bd{P}_n) Z^{-1}_{n}\Bigr)
  - \frac{\I}{2} b_n \tr \Bigl( \partial_{\bd{z}}\bd{Q}_{n}
  \nabla^2 U(\bd{Q}_n) Z^{-1}_{n} \Bigr),
\end{equation}
with initial condition $b_n(0, \bd{q}, \bd{p}) = 2^{d/2}$. The
evolution equation \eqref{eq:bamplitude} for $b_n$ is \emph{manifestly gauge-invariant}, as all terms are independent of the gauge choice. Using the amplitude function $b_{n}$, the frozen Gaussian
approximation can be rewritten as
\begin{multline}\label{eq:FGAapprox}
  \psi^{\veps}_{\FGA}(t,\bd{x})=\frac{1}{(2\pi\veps)^{3d/2}}\sum_{n=1}^{\infty} \int_{\Gamma^{\ast}}
  \int_{\RR^d} b_n(t, \bd{q}, \bd{p}) u_n\left(\bd{P}_n,{\bd{x}}/{\veps}\right) G_{\bd{Q}_n, \bd{P}_n}^{\veps}(\bd{x})
  e^{\I S_n(t, \bd{q}, \bd{p})/\veps - \I S_n^{\mc{A}}(t, \bd{q}, \bd{p})} \\ \times \average{ G_{\bd{q}, \bd{p}}^{\veps} u_n(\bd{p}, \cdot / \veps)| \psi_0} \ud \bd{q} \ud \bd{p}.
\end{multline}
The gauge-dependent term in \eqref{eq:FGAapprox} thus reads
\begin{equation}\label{eq:gaugeterms}
  u_n(\bd{P}_n, \bd{x}/\veps) e^{- \I S_{n}^{\bd{A}}(t, \bd{q}, \bd{p})}
  \overline{u}_n(\bd{p}, \bd{y}/\veps).
\end{equation}
Our goal is hence to design a gauge-invariant time
integrator for \eqref{eq:defSA} such that the term
\eqref{eq:gaugeterms} becomes independent of the gauge. Observe that,
by the Hamiltonian flow \eqref{eq:Hamiltonian Flow},
\begin{equation}\label{eq:SA}
  S_n^{\mc{A}}(t, \bd{q}, \bd{p}) = -
  \int_0^t \mc{A}(\bd{P}_n) \cdot \ud \bd{P}_n(s).
\end{equation}
Let $0 = t_0 < t_1 < \cdots < t_K = t$ be a time
discretization, we have
\begin{equation}\label{eq:expapp}
  \exp(-\I S_n^{\mc{A}}) = \exp\biggl( \I \int_0^t \mc{A}(\bd{P}_n) \cdot \ud
  \bd{P}_n(s)\biggr) = \prod_{k = 1}^{K} \exp\biggl( \I \int_{t_{k-1}}^{t_k}
  \mc{A}(\bd{P}_n) \cdot \ud \bd{P}_n(s) \biggr).
\end{equation}
To proceed, let us first work in a gauge where $u_{n}(\bd{\xi},\cdot)$
is smooth in $\bd{\xi}\in\Gamma^{\ast}$. Note that since our final
formula is gauge-independent, the choice of the gauge here is only for
the derivation. Using the Taylor approximation, we obtain
\begin{equation}\label{eq:expapprox}
  \begin{aligned}
    \I \int_{t_{k-1}}^{t_k} \mc{A}(\bd{P}_n) \cdot \ud
    \bd{P}_n(s) & = -\I\,\Im\left\{\langle u_{n}(\bd{P}_{n}(t_{k-1}),\cdot)|\nabla u_{n}(\bd{P}_{n}(t_{k-1}),\cdot)\rangle\cdot\Delta\bd{P}_{k,n}\right\}+\mathcal{O}(\Delta \bd{P}_{k,n})^2\\
    & = \I\,\Im\left\{ 1-\langle u_{n}(\bd{P}_{n}(t_{k-1}),\cdot)|u_{n}(\bd{P}_{n}(t_{k}),\cdot)\rangle\right\}+\mathcal{O}(\Delta \bd{P}_{k,n})^2\\
    & = \I \,\Im \{ \ln \bigl\langle u_{n}(\bd{P}_n(t_k), \cdot)|
    u_{n}(\bd{P}_n(t_{k-1}), \cdot)
    \bigr\rangle\}+\mathcal{O}(\Delta \bd{P}_{k,n})^2, \\
  \end{aligned}
\end{equation}
where $\Delta\bd{P}_{k,n}=\bd{P}_{n}(t_{k})-\bd{P}_{n}(t_{k-1})$. The
first approximation was obtained by using a left Riemann sum. The next
approximation is the forward difference approximation for the
derivative. The last approximation is the Taylor series for
$\ln \bd{z}$ around $\bd{z}=1$. Therefore, taking exponential, we get
\begin{equation}
\begin{aligned}
\exp\Bigl(\I \int_{t_{k-1}}^{t_k} \mc{A}(\bd{P}_n) \cdot \ud
    \bd{P}_n(s) \Bigr) & = \frac{ \bigl\langle
    u_{n}(\bd{P}_{n}(t_k), \cdot)| u_{n}(\bd{P}_{n}(t_{k-1}), \cdot)
    \bigr\rangle}{\bigl\lvert \bigl\langle
    u_{n}(\bd{P}_n(t_k), \cdot)| u_{n}(\bd{P}_n(t_{k-1}), \cdot)
    \bigr\rangle \bigr\rvert}+\mathcal{O}(\Delta \bd{P}_{k,n})^2.\\
\end{aligned}
\end{equation}
Substituting the last equation in the right hand side of
\eqref{eq:expapp} gives an approximation to
$\exp(-\I S_{n}^{\mathcal{A}})$ with and error
$\mathcal{O}(\Delta \bd{P}_{n})$ with
$\Delta \bd{P}_{n}=\displaystyle\max_{k}|\Delta \bd{P}_{k,n}|$.  
This then gives the approximation to \eqref{eq:gaugeterms} as
\begin{multline}\label{eq:Fn}
    u_n(\bd{P}_n, \bd{x}/\veps) e^{- \I S_{n}^{\bd{A}}(t, \bd{q}, \bd{p})}
  \overline{u}_n(\bd{p}, \bd{y}/\veps)
  = F_{n}(t,\bd{q},\bd{p},\bd{x},\bd{y})+ \mathcal{O}(\Delta \bd{P}_{n}) \\
  :=  \bigl \lvert u_{n}(\bd{P}_n(t_K), \bd{x}/\veps) \bigr\rangle
  \prod_{k=1}^{K}
  \frac{ \bigl\langle
    u_{n}(\bd{P}_n(t_k), \cdot)\big| u_{n}(\bd{P}_n(t_{k-1}), \cdot)
    \bigr\rangle}{\bigl\lvert \bigl\langle
    u_{n}(\bd{P}_n(t_k), \cdot)\big| u_{n}(\bd{P}_n(t_{k-1}), \cdot)
    \bigr\rangle \bigr\rvert}
  \bigl \langle u_{n}(\bd{P}_n(t_0), \bd{y}/\veps) \bigr\rvert+\mathcal{O}(\Delta \bd{P}_{n}).
\end{multline}
The right hand side of \eqref{eq:Fn} is manifestly gauge-invariant, as
the phase term in $\ket{u_{n}(\bd{P}_{n}(t_{k}),\cdot)}$ will cancel
with that of $\bra{u_{n}(\bd{P}_{n}(t_{k}),\cdot)}$, for
$k = 0, \ldots, K$.

Therefore, in summary, we arrive at a gauge-invariant reformulation of
$\psi_{\FGA}^{\veps}$ as
\begin{multline}\label{eq:GaugeFreeFGA}
\psi_{\FGA}^{\veps}(t,\bd{x})\approx\dfrac{1}{(2\pi\veps)^{3d/2}}\displaystyle\sum_{n=1}^{\infty}\displaystyle\int_{\Gamma^{\ast}}
\displaystyle\int_{\mathbb{R}^{d}}b_{n}(t,\bd{q},\bd{p})F_{n}(t,\bd{q},\bd{p},\bd{x},\bd{y})G_{\bd{Q}_{n},\bd{P}_{n}}^{\veps}(\bd{x}) \\
\times e^{\I S_{n}(t,\bd{q},\bd{p})/\veps}\braket{G_{\bd{q},\bd{p}}^{\veps}\big|\psi_{0}}\ud\bd{q}\ud\bd{p},
\end{multline}
where $F_n$ is given by \eqref{eq:Fn}, and the evolution of $(\bd{Q}_{n},\bd{P}_{n})$ follows the Hamiltonian dynamics
\begin{equation}\label{eq:Flow1}
  \left\{
  \begin{aligned}
    &\dfrac{d\bd{Q}_{n}}{dt}=\nabla E_{n}(\bd{P}_{n}),\\
    &\dfrac{d\bd{P}_{n}}{dt}=-\nabla U(\bd{Q}_{n}),
  \end{aligned}
\right.
\end{equation}
with initial condition $\bd{Q}_{n}(0,\bd{q},\bd{p})=\bd{q}$ and $\bd{P}_{n}(0,\bd{q},\bd{p})=\bd{p}$.

The action $S_n$ solves
\begin{equation}\label{eq:Flow2}
\dfrac{dS_{n}}{dt}=\bd{P}_{n}\cdot\nabla_{\bd{P}}h_{n}(\bd{Q}_{n},\bd{P}_{n})-h_{n}(\bd{Q}_{n},\bd{P}_{n}),
\end{equation}
with initial condition $S_{n}(0,\bd{q},\bd{p})=0$, and the amplitude $b_n$ follows the evolution
\begin{equation}\label{eq:Flow3}
\dfrac{db_{n}}{dt}=\dfrac{1}{2}b_{n}\text{tr}\left(\partial_{z}\bd{P}_{n}\nabla^{2}E_{n}(\bd{P}_{n})Z_{n}^{-1}\right)-\dfrac{\I}{2}b_{n}
\text{tr}\left(\partial_{z}\bd{Q}_{n}\nabla^{2}U_{n}(\bd{Q}_{n})Z_{n}^{-1}\right),
\end{equation}
with initial condition $b_{n}(0,\bd{q},\bd{p})=2^{d/2}$.

\subsection{Analytical Results}\label{results}

To make the presentation self-contained, we briefly recall here the
analytical results proved in \cite{DeLuYa:analysis} for the frozen Gaussian
approximation to \eqref{eq:schrodinger}. The proofs of these results
and more details can be found in \cite{DeLuYa:analysis}.

First we recall that the FGA ansatz recover the initial condition at
time $t=0$, $\psi_{\FGA}^{\veps}(0,\bd{x}) = \psi_0(\bd{x})$. This
follows from the Bloch decomposition \eqref{eq:Band Decomposition}.

Let us recall a few notions from \cite{DeLuYa:analysis} to state the
convergence results for the frozen Gaussian approximation. We define
the windowed Bloch transform
$(\mathcal{W}f)_{n}(\bd{q},\bd{p}):L^{2}(\mathbb{R}^{d})\rightarrow
L^{2}(\Omega)^{\mathbb{N}}$ as
\begin{equation}
(\mathcal{W}{f})_{n}(\bd{q},\bd{p}) =  \dfrac{2^{d/4}}{(2\pi)^{3d/4}}\langle u_{n}(\bd{p},\cdot)G_{\bd{q},\bd{p}}|f\rangle=\dfrac{2^{d/4}}{(2\pi)^{3d/4}}\int_{\RR^d} \wb{u}_n(\bd{p},\bd{x})\wb{G}_{\bd{q},\bd{p}}(\bd{x})
  f(\bd{x}) \ud\bd{x}
\end{equation}
where
\begin{equation}
G_{\bd{q},\bd{p}}(\bd{x}):=\exp\left(-\dfrac{|\bd{x}-\bd{q}|^{2}}{2}+\I\bd{p}\cdot(\bd{x}-\bd{q})\right).
\end{equation}
The adjoint operator $\mathcal{W}^{\ast}:L^{2}(\Omega)^{\mathbb{N}}\rightarrow L^{2}(\mathbb{R}^{d})$ is then
\begin{equation}
(\mathcal{W}^{\ast}g)(\bd{x})=\dfrac{2^{d/4}}{(2\pi)^{3d/4}}\sum_{n=1}^{\infty}\iint_{\Omega}u_{n}(\bd{p},\bd{x})G_{\bd{q},\bd{p}}(\bd{x})g_{n}(\bd{q},\bd{p})\ud\bd{q}\ud\bd{p}.
\end{equation}

The windowed Bloch transform and its adjoint have the following
important property.
\begin{prop*}[\cite{DeLuYa:analysis}*{Proposition 2.2}] The windowed Bloch transform and its adjoint satisfies
  \begin{equation}\label{eq:reconstruction}
    \mc{W}^{\ast} \mc{W} = \mathrm{Id}_{L^2(\RR^d)}.
  \end{equation}
\end{prop*}

\begin{remark}
  Similar to the windowed Fourier transform, the representation given
  by the windowed Bloch transform is redundant, so that $\mc{W}
  \mc{W}^{\ast} \not= \mathrm{Id}_{L^2(\Omega)^{\NN}}$. The
  normalization constant in the definition of $\mc{W}$ is also due to
  this redundancy.
\end{remark}

\medskip
The previous proposition motivates us to consider the contribution of
each band to the reconstruction formulae
\eqref{eq:reconstruction}. This gives to the operator $\Pi_n^{\mc{W}}:
L^2(\RR^d) \to L^2(\RR^d)$ for each $n \in \NN$
\begin{equation}
  (\Pi_n^{\mc{W}} f)(\bd{x}) = \frac{2^{d/4}}{(2\pi)^{3d/4}}
  \iint_{\Omega} u_n(\bd{p}, \bd{x}) G_{\bd{q},\bd{p}}(\bd{x})
  (\mc{W} f)_{n}(\bd{q},\bd{p})  \ud\bd{q}\ud\bd{p}.
\end{equation}
It follows from \eqref{eq:reconstruction} that $\sum_n \Pi_n^{\mc{W}} = \mathrm{Id}_{L^2(\RR^d)}$.

Correspondingly, the semiclassical windowed Bloch transform $\mc{W}^{\veps}:
L^2(\RR^d) \to L^2(\Omega)^{\NN}$ is defined as
\begin{equation}
  (\mc{W}^{\veps} f)_{n}(\bd{q},\bd{p}) =  \frac{2^{d/4}}{(2\pi\veps)^{3d/4}} \braket{u_n(\bd{p}, \cdot / \veps) G^{\veps}_{\bd{q}, \bd{p}}| f} = \frac{2^{d/4}}{(2\pi\veps)^{3d/4}} \int_{\RR^d} \wb{u}_n(\bd{p}, \bd{x}/\veps)  \wb{G}^{\veps}_{\bd{q},\bd{p}}(\bd{x}) f(\bd{x}) \ud \bd{x}.
\end{equation}
Similarly we also have the  operator $\Pi_n^{\mc{W}, \veps}:
L^2(\RR^d) \to L^2(\RR^d)$ for each $n \in \NN$ with semiclassical scaling
\begin{equation}\label{eq:projection}
  (\Pi_n^{\mc{W}, \veps} f)(\bd{y}) = \frac{2^{d/4}}{(2\pi\veps)^{3d/4}}
  \iint_{\Omega} u_n(\bd{\xi}, \bd{y}/\veps) G^{\veps}_{\bd{x},\bd{\xi}}(\bd{y})
  (\mc{W}^{\veps} f)_{n}(\bd{x},\bd{\xi})  \ud\bd{x}\ud\bd{\xi}.
\end{equation}
It follows from \eqref{eq:reconstruction} and a change of variable
that $\sum_n \Pi_n^{\mc{W}, \veps} = \mathrm{Id}_{L^2(\RR^d)}$.

For the long time existence of the Hamiltonian flow
\eqref{eq:Hamiltonian Flow}, we will assume that the external
potential $U(\bd{x})$ is subquadratic, such that
$\norm{\partial_{\bd{x}}^{\alpha} U(\bd{x})}_{L^{\infty}}$ is finite
for all multi-index $|\alpha|\geq 2$.  As a result, since the domain
of $\bd{\xi}$ is bounded, the Hamiltonian $h_n$ is also
subquadratic. $\psi_{FGA}^{\veps}$ provides an approximate solution to
equation \eqref{eq:schrodinger} to first order accuracy as stated in the
two theorems below, rephrased from our previous work \cite{DeLuYa:analysis}.

\begin{theorem*}[\cite{DeLuYa:analysis}*{Theorem 3.1}]
  Assume that the $n$-th Bloch band $E_{n}(\bd{\xi})$ does not
  intersect any other Bloch bands for all $\bd{\xi}\in \Gamma^{\ast}$;
  and moreover, the Hamiltonian $h_n(\bd{x}, \bd{\xi})$ is
  subquadratic. Let $\mathscr{U}^{\veps}_t$ be the propagator of
  the time-dependent Schr\"{o}dinger equation \eqref{eq:schrodinger}. Then for any given $T$, $0 \leq t \leq T$ and sufficiently small,
  $\veps\leq \veps_0$,
  \begin{multline}\label{eq:theoremA}
    \sup_{0\leq t\leq T} \, \Bigl\lVert\, \mathscr{U}^{\veps}_t \bigl(
    \Pi^{\mc{W}, \veps}_n \psi_0\bigr) - \frac{1}{(2\pi\veps)^{3d/2}}
    \int_{\Omega} b_n(t, \bd{q}, \bd{p}) u_n(\bd{P}_n, \bd{x}/ \veps)
    G_{\bd{Q}_n, \bd{P}_n}^{\veps}(\bd{x}) \times \\
    \times e^{\I S_n(t, \bd{q}, \bd{p})/\veps - \I S_n^{\mc{A}}(t,
      \bd{q}, \bd{p})} \average{ G_{\bd{q}, \bd{p}}^{\veps}
      u_n(\bd{p}, \cdot / \veps)| \psi_0} \ud \bd{q} \ud \bd{p} \,
    \Bigr\rVert_{L^{2}} \leq C_{T, n} \, \veps\, \bigl\lVert
    \psi_{0}^{\veps} \bigr\rVert_{L^2}.
  \end{multline}
\end{theorem*}
\begin{theorem*}[\cite{DeLuYa:analysis}*{Theorem 3.2}]
  Assume that the first $N$ Bloch bands $E_{n}(\bd{\xi})$,
  $n=1,\cdots,N$ do not intersect and are separated from the other
  bands for all $\bd{\xi}\in \Gamma^{\ast}$; and assume that the
  Hamiltonian $h_n(\bd{x}, \bd{\xi})$ is subquadratic. Let
  $\mathscr{U}^{\veps}_t$ be the propagator of the time-dependent
  Schr\"{o}dinger equation~\eqref{eq:schrodinger}. Then for any given
  $T$, $0 \leq t \leq T$ and sufficiently small $\veps$, we have
  \begin{multline}
    \sup_{0\leq t\leq T} \, \biggl\lVert\,\mathscr{U}^{\veps}_t
    \psi_0 -
    \frac{1}{(2\pi\veps)^{3d/2}}  \sum_{n=1}^{N}
    \int_{\Omega} b_n(t, \bd{q}, \bd{p}) u_n(\bd{P}_n, \bd{x}/ \veps)
    G_{\bd{Q}_n, \bd{P}_n}^{\veps}(\bd{x}) e^{\I S_n(t, \bd{q}, \bd{p})/\veps - \I S_n^{\mc{A}}(t,
      \bd{q}, \bd{p})}  \times \\
    \times \average{ G_{\bd{q}, \bd{p}}^{\veps}
      u_n(\bd{p}, \cdot / \veps)| \psi_0} \ud \bd{q} \ud \bd{p}  \, \biggr\rVert_{L^{2}}
    \leq C_{T, N}\, \veps\bigl\lVert
    \psi_{0}^{\veps}\bigr\rVert_{L^2}+ \norm{ \psi_0^{\veps} -
      \sum_{n=1}^N \Pi^{\mc{W}, \veps}_n \psi_0^{\veps}}_{L^2}.
  \end{multline}
\end{theorem*}

These approximation results show the first order asymptotic accuracy of FGA, which
will be numerically validated in  Section~\ref{sec:example}.

\section{Numerical implementation}\label{sec:algorithm}

We will now describe the numerical implementation of the
gauge-invariant frozen Gaussian approximation (GIFGA) method.
We will restrict ourselves to one spatial dimension in this paper.
For one thing, the computation of true solutions to \eqref{eq:schrodinger}
with high accuracy is extremely time-consuming in high dimensions,
and thus it is difficult for us to confirm numerically the asymptotic
convergence order with the pollution of non-negligible numerical errors.
For another thing, band-crossing is quite common in high dimensional cases
(e.g., in honeycomb lattice), which requires more techniques than
the scope of this paper, and we will leave the numerical study of high dimensional
examples as future work. The calculation of the
Bloch eigenvalues and eigenfunctions is discussed in
Section~\ref{blochbands}. 
In Section~\ref{algorithmdescription} We describe the numerical
algorithms of GIFGA based on the Bloch bands. We will also
discuss the mesh sizes required for accurate computation.

\subsection{Numerical computation of Bloch bands}\label{blochbands}
We show how to compute numerically the eigenvalues and eigenfunctions of \eqref{eq:eigen} in $d=1$.
Define the Fourier transform of $u_{n}(\xi,x)$ as
\begin{equation}
  \widehat{u}_{n}(\xi,\eta)=\dfrac{1}{2\pi}\int_{\Gamma}u_{n}(\xi,x)e^{-\I x\eta}\ud x.
\end{equation}
Taking the Fourier transform of \eqref{eq:eigen} one obtains
\begin{equation}
\dfrac{(\eta+\xi)^{2}}{2}\wh{u}_{n}(\xi,\eta)+\wh{V_\Gamma}(\eta)*\wh{u_n}(\xi,\eta)=E_{n}(\xi)\wh{u}_{n}(\xi,\eta),
\end{equation}
where ``$*$'' stands for the operation of convolution.

Truncating the Fourier grid to
$\left\{-\Lambda,\cdots,\Lambda-1\right\}\subset\mathbb{Z}$ gives
\begin{equation}
  H_{\xi}(\Lambda)\left(\begin{matrix}
      \wh{u}_{n}(\xi,-\Lambda)\\
      \wh{u}_{n}(\xi,1-\Lambda)\\
      \vdots\\
      \wh{u}_{n}(\xi,\Lambda-1)\\
\end{matrix}\right)=E_{n}(\xi)\left(\begin{matrix}
\wh{u}_{n}(\xi,-\Lambda)\\
\wh{u}_{n}(\xi,1-\Lambda)\\
\vdots\\
\wh{u}_{n}(\xi,\Lambda-1)\\
\end{matrix}\right)
\end{equation}
where $H_{\xi}(\Lambda)$ is the $2\Lambda\times 2\Lambda$ matrix given by
\begin{equation}
  H_{\xi}(\Lambda)=\left(\begin{matrix}
      \dfrac{(-\Lambda+\xi)^{2}}{2}+\wh{V}_{\Gamma}(0) & \wh{V}_{\Gamma}(-1) & \cdots & \wh{V}_{\Gamma}(1-2\Lambda)\\
      \wh{V}_{\Gamma}(1) & \dfrac{(-\Lambda+1+\xi)^{2}}{2}+\wh{V}_{\Gamma}(0) & \cdots & \wh{V}_{\Gamma}(2-2\Lambda)\\
      \hdots & \hdots & \ddots & \hdots\\
      \wh{V}_{\Gamma}(2\Lambda-1) & \wh{V}(2\Lambda-2) & \cdots & \dfrac{(\Lambda-1+\xi)^{2}}{2}+\wh{V}_{\Gamma}(0)\\
\end{matrix}\right).
\end{equation}
After diagonalizing the matrix, the eigenfunction in the physical domain is then obtained via inverse Fourier transform
\begin{equation}
  u_{n}(\xi, x)\approx\sum_{y=-\Lambda}^{\Lambda-1}\wh{u}_{n}(\xi,\eta)e^{\I \eta x}.
\end{equation}

\begin{figure}
\includegraphics[scale=.6]{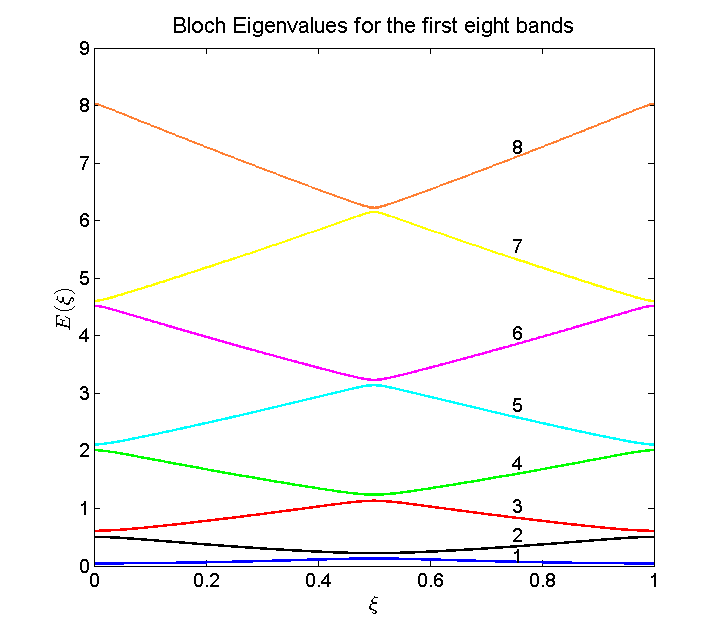}
\caption{Energy eigenvalues for the one-dimensional lattice potential $V(x)=\exp\left(-25x^{2}\right)$}\label{fig1}
\end{figure}

\begin{figure}
\includegraphics[scale=.6]{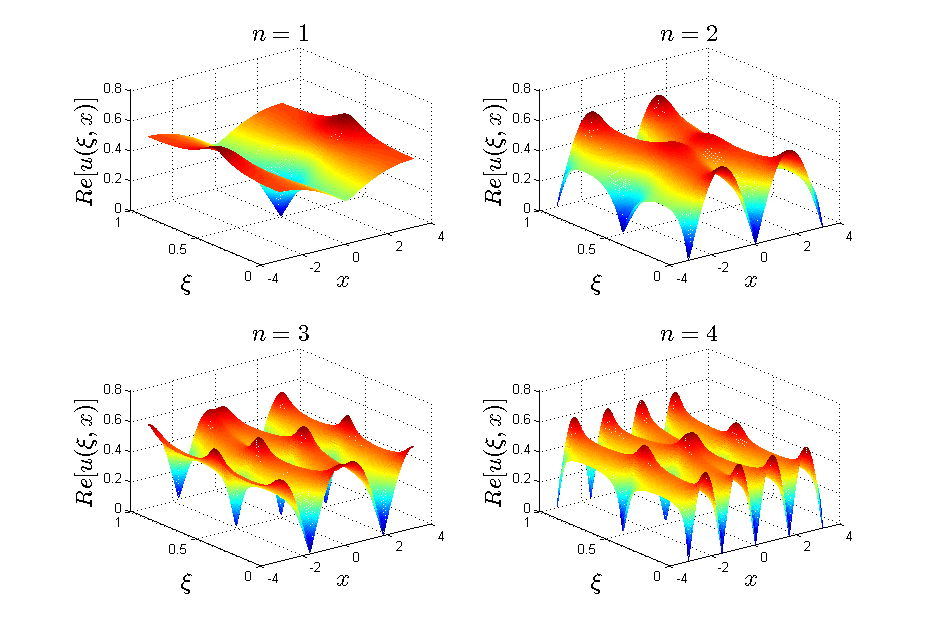}
\caption{Module of eigenfunctions for the one-dimensional lattice potential $V(x)=\exp\left(-25x^{2}\right)$. We display absolute value of the first 4 lowest energy eigenfunctions.}\label{fig2}
\end{figure}

\begin{figure}
\hspace{5em}\includegraphics[scale=.65]{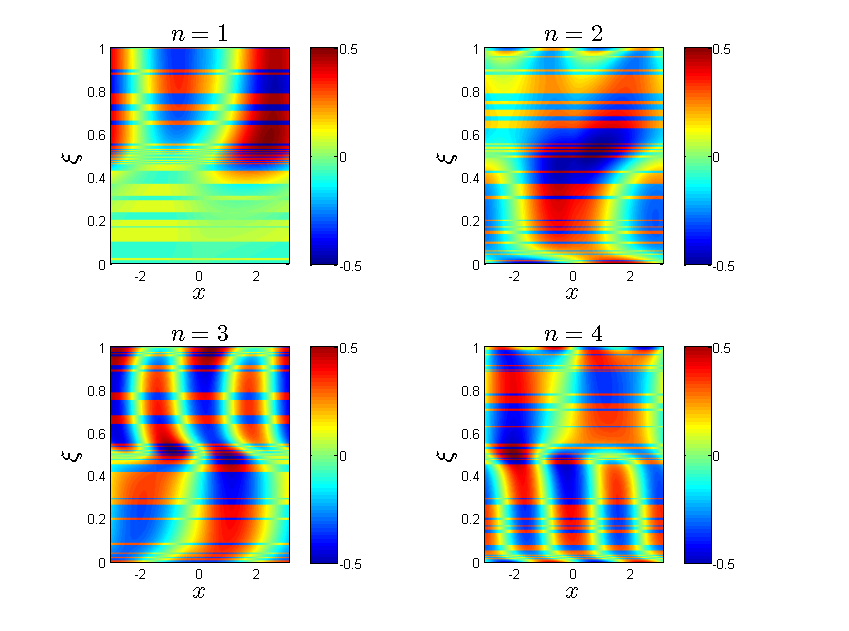}
\caption{Real part of the eigenfunctions for the one-dimensional lattice potential $V(x)=\exp\left(-25x^{2}\right)$. We display the real parts for the first 4 lowest energy eigenfunctions. We use 100 data points for the $\bd{\xi}$ variable.}\label{fig3}
\end{figure}

\begin{figure}
\hspace{5em}\includegraphics[scale=.65]{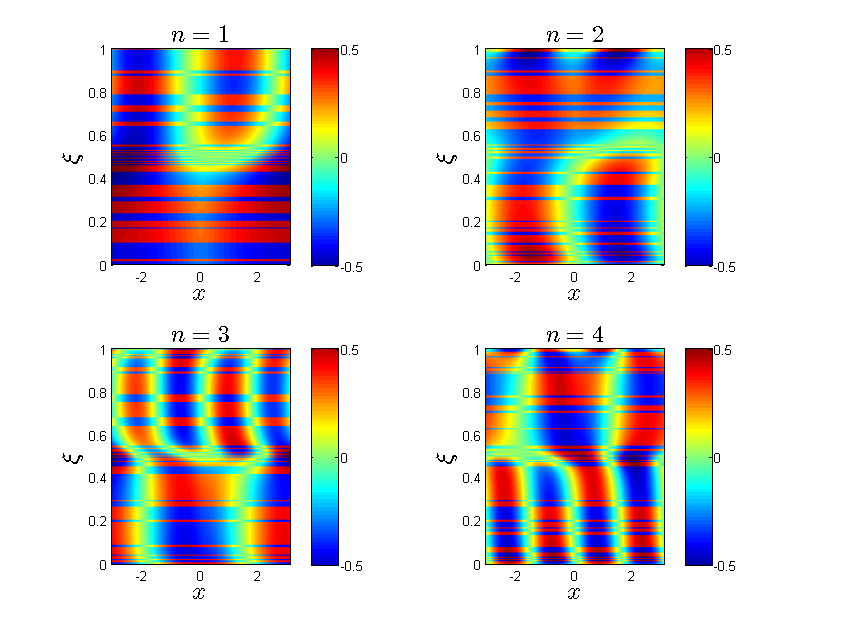}
\caption{Imaginary part of the eigenfunctions for the one-dimensional lattice potential $V(x)=\exp\left(-25x^{2}\right)$. We display the imaginary parts for the first 4 lowest energy eigenfunctions. We use 100 data points for the $\bd{\xi}$ variable.}\label{fig4}
\end{figure}

\begin{example}\label{ex1}
  In this example, we compute Bloch eigenvalues and eigenfunctions
  with potential $V_{\Gamma}(x)=\exp\left(-25x^{2}\right)$. The
  extension of $V_{\Gamma}(x)$ periodically with respect to $\Gamma$
  is not analytic on the boundary of $\Gamma$. However, this lack of
  smoothness presents a negligible problem numerically as
  $V_{\Gamma}(x)$ decays rapidly.  Figure~\ref{fig1} shows the energy
  eigenvalues $E_{n}(\xi)$ for $\xi\in [0,1)$. The plot shows the
  first $8$ bands where the bottom curve corresponds to $n=1$ (lowest
  band) and the top curve represents $n=8$ (highest band). Figure
  \ref{fig2} shows the modules of the corresponding Bloch
  eigenfunctions for the first $4$ bands. Notice that while these
  surfaces are continuous and periodic, the next two figures
  (\ref{fig3} and \ref{fig4}) of the real and imaginary parts of the
  Bloch eigenfunctions are not. This is due to the arbitrary gauge
  freedom in the diagonalization.
\end{example}

\remark 1. In the numerical computation of $E(\xi)$, the corresponding
eigenfunctions and their derivatives near the points $\xi=0$ and
$\xi=0.5$ (and $\xi=1$ by periodicity) is tricky, since the Bloch
bands are close to each other near these points (see Figure~\ref{fig1}). For this reason, our grid for the $\xi$
variable will not contain these points. In other words, we shift the
grids in the first Brillouin zone to avoid these high symmetry points.

2. One can apply the same technique to derive an algorithm for computing Bloch eigenvalues and eigenfunctions in higher dimensions. The main issue with this algorithm is that the numerical cost increases drastically for $d>1$. In the case where the periodic potential has the form $V_{\Gamma}(\bd{x})=\sum_{j=1}^{d}V_{j}(x_{j})$ with $V_{j}(x_{j}+2\pi)=V(x_{j})$, computation of Bloch bands can be treated for each coordinate $x_{j}$ separately. For some common potentials, data for the energy eigenvalues has already been produced (see remark 2.1 in \cite{HuJiMaSp:07}).



\subsection{Algorithms for gauge invariant frozen Gaussian approximation}\label{algorithmdescription}


We assume that the initial data $\psi_{0}(x)$ has compact support or
that it decays sufficiently fast as $\abs{x}\rightarrow\infty$, and
hence, we only need to use a finite number of mesh points in physical
space.

For a mesh size $\delta x$ and starting point
 $x^{0} \in \RR$, the grid is specified as
\begin{equation}
  x^{m}= x^{0}+(m-1)\delta x,
\end{equation}
for $m=1,\cdots,N_{x}$, where $N_x$ is the number of the spatial grid
in one dimension.

We present the algorithm in five steps below.

\noindent \textbf{Step 1}. Compute the Bloch eigenvalues $E_{n}(\xi)$ and eigenfunctions $u_{n}(\xi,x)$ of \eqref{eq:eigen}, according to the
algorithm described in Section~\ref{blochbands}.

\remark For our one dimensional examples in Section \ref{sec:example},
we choose a mesh for $(\xi,x)$ such that $\delta\xi=(1-2\rho)/199$
with $\xi^{0}=-1/2+\rho$ and $N_{\xi}=200$; and $\delta x=2\pi/804$
with $x^{0}=-\pi$ and $N_{x}=805$ for some $0<\rho\ll 1$. $\rho$ was
included to avoid putting mesh points at high symmetry points in the
first Brillouin zone. This number of grid points is enough to ensure
that the eigenvalues and eigenfunctions are computed with sufficient
accuracy for our numerical tests.

\smallskip

\noindent \textbf{Step 2}. Compute $(Q_{n}(t,q,p),
P_{n}(t,q,p),S_{n}(t,q,p),b_{n}(t,q,p))$ in \eqref{eq:Flow1},
\eqref{eq:Flow2}, and \eqref{eq:Flow3}.

To integrate the ODEs for $(Q_{n},P_{n},S_{n},b_{n})$, we use a
symplectic fourth order Runge-Kutta method. Coefficients for the
Butcher tableau can be found in \cite{Ge00}. We will choose a mesh for
$({q},{p})\in\Omega$ and $\left({Q}_{n},{P}_{n}\right)$ takes initial
value at the grid points. That is,
\begin{align}
{Q}_{n}(0,{q},{p})={q}^{{I}}=&{q}^{0}+{I} \delta{q}\\
{P}_{n}(0,{q},{p})={p}^{{J}}=&{p}^{0}+{J} \delta{p}
\end{align}
where $I \in 1,\cdots, N_{I}$ and $J \in 1,\cdots, N_{J}$. Notice that
to represent the initial condition $\psi_{\FGA}^{\veps}(0, {x})$ one
only needs the mesh points ${q}^{{I}}$ near ${x}$. To be more precise,
as the standard deviation of the semiclassical Gaussians in
\eqref{eq:Gaussian} is $\sqrt{\veps}$ so one only needs the mesh points
$q^I$ contributing significantly to $\psi_{\FGA}^{\veps}(0, x)$
satisfy $|{x}-{q}^{{I}}|\leq\mathcal{O}(\sqrt{\veps})$. This implies
that one can put a finite number of mesh points for ${q}$-coordinate
and not on all of $\mathbb{R}$. The mesh size for ${q}^{{I}}$ and
${p}^{{J}}$ is chosen to be $\mathcal{O}(\sqrt{\veps})$, which
resolves the oscillation of the initial condition.

\smallskip
\noindent \textbf{Step 3}. Compute the windowed Bloch transformation
of the initial condition
$\braket{u_{n}({p},{\cdot}/\veps)G_{{q},{p}}^{\veps} | \psi_{0}}$. For
the sake of convenience, denote this term by
$w_{n}^{\veps}({q},{p})$. Let
\begin{equation}
{y}^{K}=y^{0}+(K-1)\delta y
\end{equation}
be a discrete mesh of ${y}$ where $K=1,\cdots,N_{y}$. Then,
\begin{equation}
w_{n}^{\veps}({q}^{{I}},{p}^{{J}})\approx\sum_{K=1}^{N_{y}}
\wb{G}_{{q}^{{I}},{p}^{{J}}}^{\veps}({y}^{{K}})\wb{u}({p}^{{J}},{y}^{{K}}/\veps)
\psi_{0}({y}^{{K}})r_{\theta}\left(|{y}^{{K}}-{q}^{{I}}|\right)\delta{{y}},
\end{equation}
with $r_{\theta}$  a cut-off function such that $r_{\theta}=1$ in the ball of radius $\theta>0$ centered at the origin and $r_{\theta}=0$ outside the ball.

The mesh ${y}^{{K}}$ should approximately cover the support of the
initial condition $\psi_{0}({y})$. As can be seen by the form of
$w_{n}^{\veps}$, the size of $N_{{y}}$ will depend on $\veps$. The
mesh should be fine enough to accurately capture
$\wb{u}_{n}({p},{y}/\veps)\wb{G}_{{q},{p}}^{\veps}({y})\psi_{0}({y})$
for all bands $n$.

\remark One can reduce the computation time of
$w_{n}^{\veps}({q}^{{I}},{p}^{{J}})$ by incorporating the periodicity
of $u_{n}({\xi},{x})$ with respect to ${x}$. As can be seen by
Figure~\ref{fig2}, $u_{n}({\xi},{x})$ tends to become more oscillatory
as $n$ increases. Thus, the mesh of ${y}^{{K}}$ should be adapted so
that it depends on $n$.

\smallskip
\noindent \textbf{Step 4}. Denote the product term in \eqref{eq:Fn} by
\begin{equation}\label{eq:F}
\wt{F}_{n}(t,{q},{p}):= \prod_{k=1}^{K}
  \frac{ \bigl\langle
    u({P}_n(t_k), \cdot), u({P}_n(t_{k-1}), \cdot)
    \bigr\rangle}{\bigl\lvert \bigl\langle
    u({P}_n(t_k), \cdot), u({P}_n(t_{k-1}), \cdot)
    \bigr\rangle \bigr\rvert},
\end{equation}
and note that
\begin{equation*}
F_n(t,q,p,x,y)=\bigl \lvert u_{n}({P}_n(t_K), {x}/\veps) \bigr\rangle
 \wt{F}_{n}(t,{q},{p})
  \bigl \langle u_{n}({P}_n(t_0), {y}/\veps) \bigr\rvert .
\end{equation*}

At this point we now have the required data to compute
$\wt{F}_{n}$. Discretize $\wt{F}_{n}$ using the same mesh from the
previous steps to obtain $\wt{F}_{n}(t,{q}^{{I}},{p}^{{J}})$. Here,
$t_{0}=0\leq t_{1}\leq t_{2}\leq\cdots\leq t_{K}=t$ is the temporal
mesh used in Step 2, with
\begin{equation*}
t_j=j\delta t,\quad \delta t=\frac{t}{K},\;\text{and}\;j=1,\cdots,K.
\end{equation*}

\smallskip
\noindent \textbf{Step 5}. Reconstruct the solution using
\eqref{eq:GaugeFreeFGA}
\begin{equation}\label{eq:discritizedFGA}
 \begin{aligned} \psi_{\FGA}^{\veps}(t,{x}^{{L}})\approx\sum_{n=1}^{N}\sum_{{I}}\sum_{{J}}&\Bigl(b_{n}(t,{q}^{{I}},{p}^{{J}})
 \wb{u}_{n}\bigl({P}_{n}(t,{q}^{{I}},{p}^{{J}}),{x}^{{L}}/\veps\bigr) G^\veps_{Q_n,P_n}(x^L)e^{S_n(t,q^I,p^J)/\veps}\\
&\times\wt{F}_{n}(t,{q}^{{I}},{p}^{{J}})\wt{\psi}_{n}^{\veps}({q}^{{I}},{p}^{{J}})r_{\theta}\left(|{x}^{{L}}-{Q}_{n}^{{I},{J}}|\right)\Bigr)\delta q\delta p,
\end{aligned}
\end{equation}
where $Q_n$ and $P_n$ are evaluated at $(t,q^I,p^J)$, and $r_{\theta}$ is a cutoff function as described in Step 3 and $N$
is the maximum number of Bloch bands used.

\noindent{\it Accuracy}. The theorems in Section~\ref{results} and \eqref{eq:Fn} imply the above algorithm has a total accuracy 
$ \mathcal{O}\left(\veps+\delta t^4/\veps+\max_n \Delta \bd{P}_{n}\right)+ \norm{ \psi_0^{\veps} -
      \sum_{n=1}^N \Pi^{\mc{W}, \veps}_n \psi_0^{\veps}}_{L^2}$, where $\mathcal{O}(\delta t^4/\veps)$ comes from the approximation 
      to the phase functions in \eqref{eq:GaugeFreeFGA} and $\mathcal{O}(\max_n \Delta \bd{P}_{n})\approx \mathcal{O}(\delta t)$ is due to the approximation \eqref{eq:Fn}. $\norm{ \psi_0^{\veps} -
      \sum_{n=1}^N \Pi^{\mc{W}, \veps}_n \psi_0^{\veps}}_{L^2}$ is the initial decomposition error, which in general decays with the number of bands as indicated in, e.g., Examples~\ref{ex2} and \ref{ex3}, and in \cite{HuJiMaSp:07}.

\section{Numerical examples}\label{sec:example}

In this section, we show the numerical performance of gauge invariant frozen Gaussian approximation (GIFGA) by several one dimensional examples, which also confirm the first order asymptotic convergency analyzed in \cite{DeLuYa:analysis}.

\subsection{Initial decomposition}
In the first two examples, we test the initial decomposition of GIFGA described in Section~\ref{sec:Bloch}. We compute $\psi_{\FGA}^{\veps}$ at $t=0$ via equation \eqref{eq:FGA}. As we cannot numerically sum to infinity, we choose to use at most $8$ bands in all of our examples. Expressed differently, the solution will be concentrated on the first 8 bands. Because of the need for $\mathcal{O}(\sqrt{\veps})$ mesh size for both coordinates $({q}^{{I}},{p}^{{J}})$ of phase space, we choose approximately $2/\sqrt{\veps}$ number of grid points for each unit interval.

\begin{example}\label{ex2}
In this example, we check the initial decomposition by choosing $\psi_{0}=A(x)\exp\bigl(\I S(x)/\veps\bigr)$ with $A(x)=\exp\left(-50x^{2}\right)\cos((x-0.5)/\veps)$ and $S(x)=0.3(x-0.5)+0.1\sin(x-0.5)$, and the lattice potential $V_{\Gamma}=\cos(x)$. We record the data in Table~\ref{ex2table}.
\end{example}

\begin{table}
\begin{tabular}{|>{\columncolor[gray]{0.8}}p{4cm}|>{\columncolor[gray]{0.8}}p{4cm}|>{\columncolor[gray]{0.8}}p{4cm}|}
\hline
$\veps=1/64$ &Error $||\psi_{0}-\psi_{\FGA}^{\veps}||_{L^2}$\\\hline
\textbf{$N=1$}&0.13260 \\\hline
\textbf{$N=2$}&0.11328\\\hline
\textbf{$N=4$}&0.033126\\\hline
\textbf{$N=8$}&7.2587e-05\\\hline
\end{tabular}

\smallskip
\begin{tabular}{|>{\columncolor[gray]{0.8}}p{4cm}|>{\columncolor[gray]{0.8}}p{4cm}|>{\columncolor[gray]{0.8}}p{4cm}|}
\hline
$\veps=1/128$ & Error $||\psi_{0}-\psi_{\FGA}^{\veps}||_{L^2}$\\\hline
\textbf{$N=1$}&0.15361 \\\hline
\textbf{$N=2$}&0.096905\\\hline
\textbf{$N=4$}&0.031652\\\hline
\textbf{$N=8$}&7.0574e-05\\\hline
\end{tabular}

\smallskip
\begin{tabular}{|>{\columncolor[gray]{0.8}}p{4cm}|>{\columncolor[gray]{0.8}}p{4cm}|>{\columncolor[gray]{0.8}}p{4cm}|}
\hline
$\veps=1/256$ &Error $||\psi_{0}-\psi_{\FGA}^{\veps}||_{L^2}$\\\hline
\textbf{$N=1$}&0.14165 \\\hline
\textbf{$N=2$}&0.1063\\\hline
\textbf{$N=4$}&0.032405\\\hline
\textbf{$N=8$}&6.9192e-05\\\hline
\end{tabular}

\smallskip
\begin{tabular}{|>{\columncolor[gray]{0.8}}p{4cm}|>{\columncolor[gray]{0.8}}p{4cm}|>{\columncolor[gray]{0.8}}p{4cm}|}
\hline
$\veps=1/512$ & Error $||\psi_{0}-\psi_{\FGA}^{\veps}||_{L^2}$\\\hline
\textbf{$N=1$}&0.15885 \\\hline
\textbf{$N=2$}&0.09276\\\hline
\textbf{$N=4$}&0.031263\\\hline
\textbf{$N=8$}&6.8701e-05\\\hline
\end{tabular}

\smallskip
\caption{$L^2$ error of $\psi_{0}({x})-\psi_{\FGA}^{\veps}(0,{x})$ for Example~\ref{ex2}. We display various values of $\veps$ and sum over $N$ Bloch bands in $\psi_{\FGA}^{\veps}$.}
\label{ex2table}
\end{table}

\begin{example}\label{ex3}
In this example, we check the initial decomposition by choosing $\psi_{0}=A(x)\exp\bigl(\I S(x)/\veps\bigr)$ with $A(x)=\exp\left(-50x^{2}\right)$ and $S(x)=0.3+0.1\sin(x-0.5)$, and the lattice potential to be $V_{\Gamma}=\exp(-25x^{2})$. We record the data in Table~\ref{ex3table}.
\end{example}

\begin{table}
\begin{tabular}{|>{\columncolor[gray]{0.8}}p{4cm}|>{\columncolor[gray]{0.8}}p{4cm}|>{\columncolor[gray]{0.8}}p{4cm}|}
\hline
$\veps=1/64$ &Error $||\psi_{0}-\psi_{\FGA}^{\veps}||_{L^2}$\\\hline
\textbf{$N=1$}&0.035736\\\hline
\textbf{$N=2$}&0.02463\\\hline
\textbf{$N=4$}&0.0075756\\\hline
\textbf{$N=8$}&0.0018796\\\hline
\end{tabular}

\smallskip
\begin{tabular}{|>{\columncolor[gray]{0.8}}p{4cm}|>{\columncolor[gray]{0.8}}p{4cm}|>{\columncolor[gray]{0.8}}p{4cm}|}
\hline
$\veps=1/128$&Error $||\psi_{0}-\psi_{\FGA}^{\veps}||_{L^2}$\\\hline
\textbf{$N=1$}&0.031445\\\hline
\textbf{$N=2$}&0.024814\\\hline
\textbf{$N=4$}&0.007579\\\hline
\textbf{$N=8$}&0.0018579\\\hline
\end{tabular}

\smallskip
\begin{tabular}{|>{\columncolor[gray]{0.8}}p{4cm}|>{\columncolor[gray]{0.8}}p{4cm}|>{\columncolor[gray]{0.8}}p{4cm}|}
\hline
$\veps=1/256$ &Error $||\psi_{0}-\psi_{\FGA}^{\veps}||_{L^2}$\\\hline
\textbf{$N=1$}&0.030633\\\hline
\textbf{$N=2$}&0.024967\\\hline
\textbf{$N=4$}&0.0076045\\\hline
\textbf{$N=8$}&0.0018698\\\hline
\end{tabular}

\smallskip
\begin{tabular}{|>{\columncolor[gray]{0.8}}p{4cm}|>{\columncolor[gray]{0.8}}p{4cm}|>{\columncolor[gray]{0.8}}p{4cm}|}
\hline
$\veps=1/512$ &Error $||\psi_{0}-\psi_{\FGA}^{\veps}||_{L^2}$\\\hline
\textbf{$N=1$}&0.030375\\\hline
\textbf{$N=2$}&0.025078\\\hline
\textbf{$N=4$}&0.0076103\\\hline
\textbf{$N=8$}&0.0018769\\\hline
\end{tabular}

\smallskip
\caption{$L^2$ error of $\psi_{0}({x})-\psi_{\FGA}^{\veps}(0,{x})$ for Example~\ref{ex3}. We display various values of $\veps$ and sum over $N$ Bloch bands in $\psi_{\FGA}^{\veps}$.}
\label{ex3table}
\end{table}

Tables \ref{ex2table}, and \ref{ex3table} show that FGA indeed matches the initial condition more closely as $N$ increases. Furthermore, we have essentially the same $L^2$ error for each $\veps$. This provides numerical verification of the independence of $\veps$ of the initial decomposition.

\remark Let us note that from equation \eqref{eq:Hamiltonian} the convergence rate should depend on the form of the lattice potential $V_{\Gamma}({x})$. Also, by equation \eqref{eq:Band Decomposition}, the convergence rate also depends on the form of the initial condition. We see from Examples \ref{ex2}, and \ref{ex3} that the cosine lattice potential seem to produce faster convergence with respect to the number of bands used. Different initial conditions may also converge faster as N increases. Example~\ref{ex5} uses an initial condition projected onto the first band. Choosing such initial condition has the advantage of needing only to compute $\psi_{\FGA}^{\veps}$ over one band.


By examining the $L^2$ errors or the convergence rates, one could determine the minimum number of bands to sum over to achieve required accuracy. In Example~\ref{ex2}, it shows that upon summing over $N=4$ bands, the initial decomposition starts to resemble the initial condition. 

\begin{figure}
\includegraphics[scale=.45]{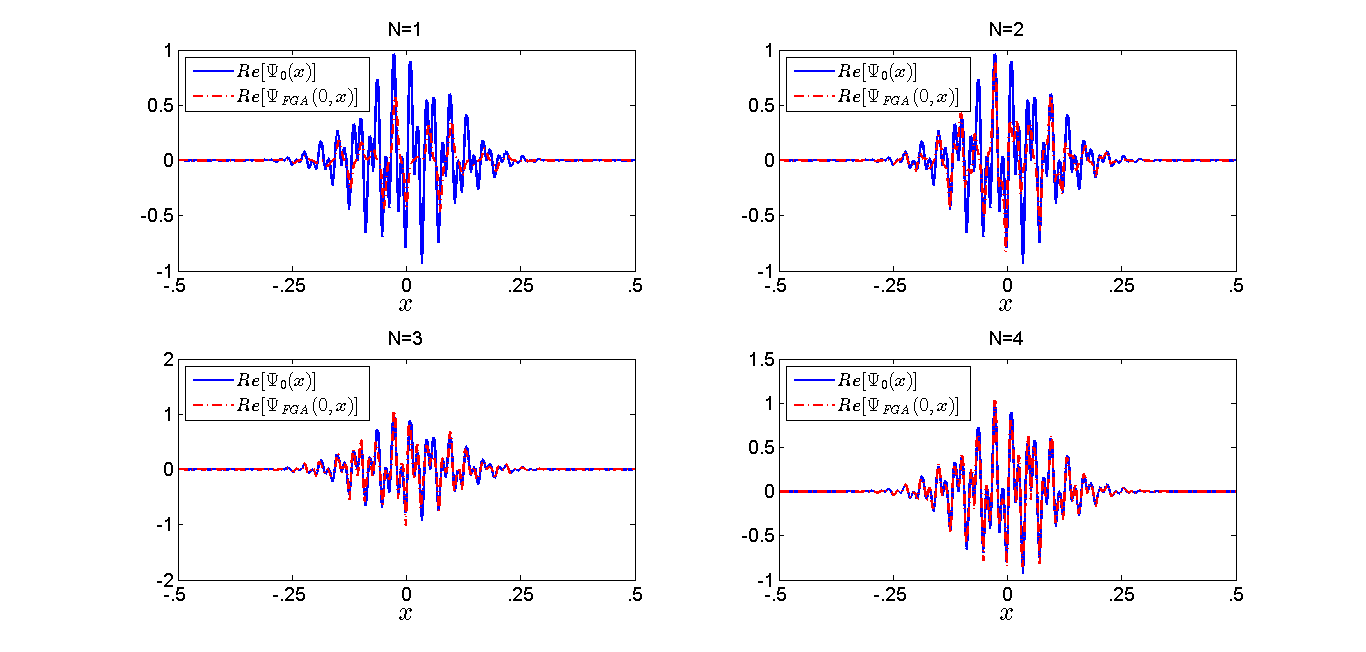}
\caption{Initial decomposition for example \ref{ex2}. The real part of $\psi_{0}({x})$ and $\psi_{\FGA}^{\veps}(0,{x})$ are shown for $\veps=1/256$. The summation in $\psi_{\FGA}^{\veps}(0,{x})$ is over the first 4 lowest energy bands.}\label{fig8}
\end{figure}

\begin{figure}
\includegraphics[scale=.45]{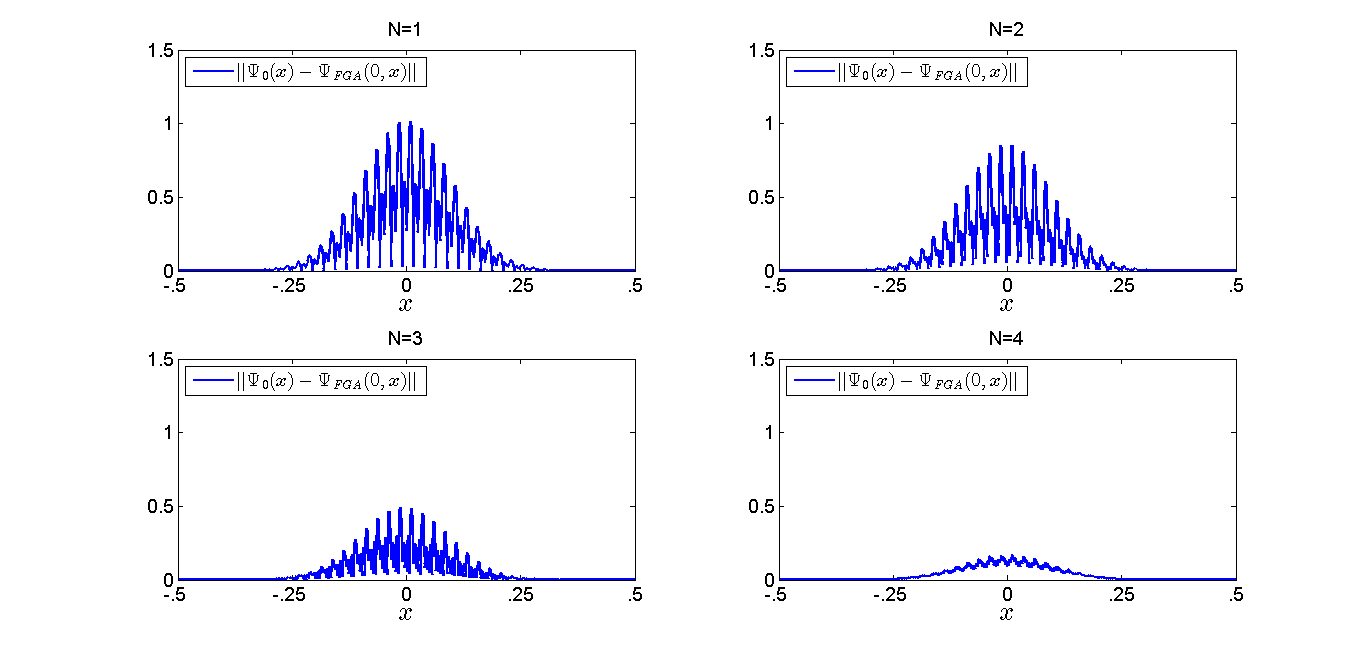}
\caption{The plot of $||\psi_{0}({x})-\psi_{\FGA}^{\veps}(0,{x})||_{l^{2}}$ for figure \ref{fig8} is displayed here.}\label{fig9}
\end{figure}

\subsection{Verification of the convergence rate of FGA}

First, we choose to test the convergence rate of \eqref{eq:FGA} with external potential $U({x})=0$ in Examples~\ref{ex4} and \ref{ex5}. With this choice of potential, there is no need for a gauge-invariant algorithm. One can optimize the algorithm described in Section~\ref{algorithmdescription} by setting $\tilde{F}(t,{q},{p})=1$ in \eqref{eq:F} in Step 4. Thus, for Examples ~\ref{ex4} and \ref{ex5}, numerical errors coming from $\tilde{F}(t,{q},{p})$ will be absent. Examples \ref{ex6} and \ref{ex7} have nonzero external potential so there will be some numerical errors introduced by $\tilde{F}(t,{q},{p})$. We continue using $2/\sqrt{\veps}$ mesh points per unit interval in ${q}$ and ${p}$ and sum up to eight bands (except for Example~\ref{ex5}). We choose a time step of size $\Delta t=T/150$. The exact solution to equation \eqref{eq:schrodinger} will be computed using the Strang splitting spectral method \cite{BaJiMa:02}. For all of our examples, the Strang splitting spectral method did not need a mesh finer than $\Delta x=1/2^{16}$ and $\Delta t=1/2^{12}$.

\begin{figure}
\includegraphics[scale=.5]{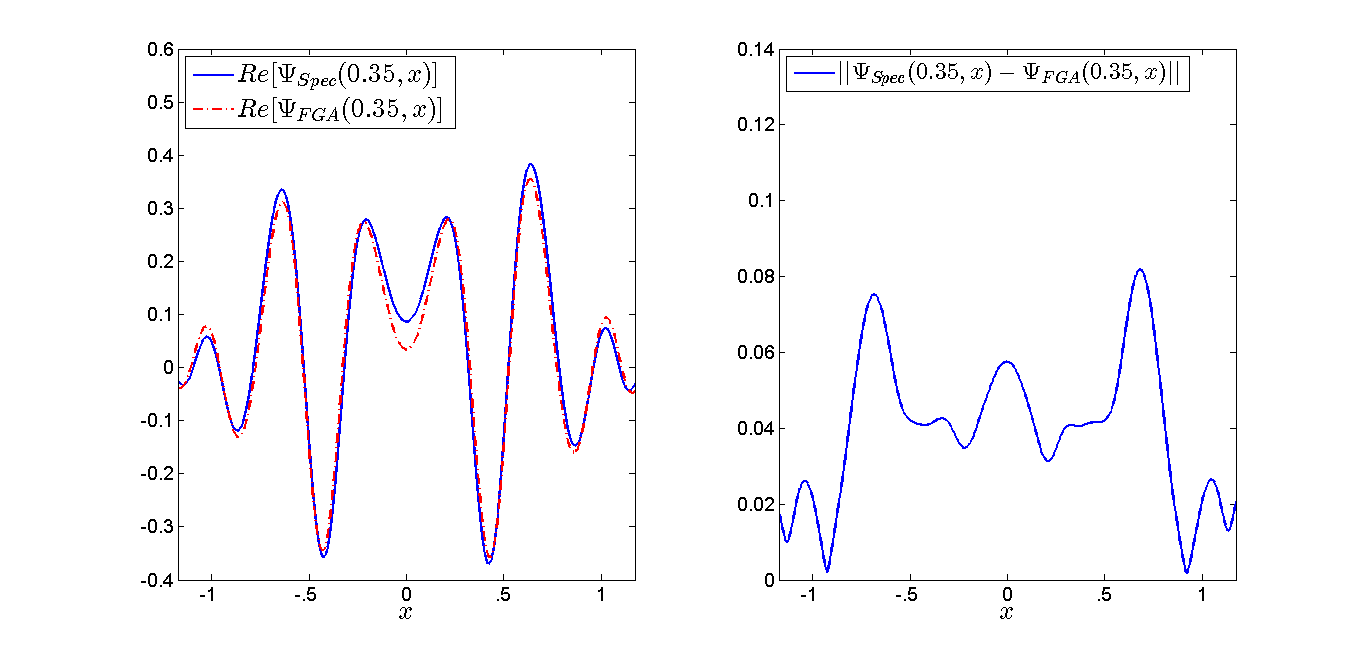}
\caption{Example \ref{ex4} plot of real parts of $\psi_{\FGA}^{\veps}(0.35,{x})$ and $\psi_{Spec}(0.35,{x})$ along side with the $L^2$ error for $\veps=1/8$.}
\label{fig16}
\end{figure}

\begin{figure}
\includegraphics[scale=.5]{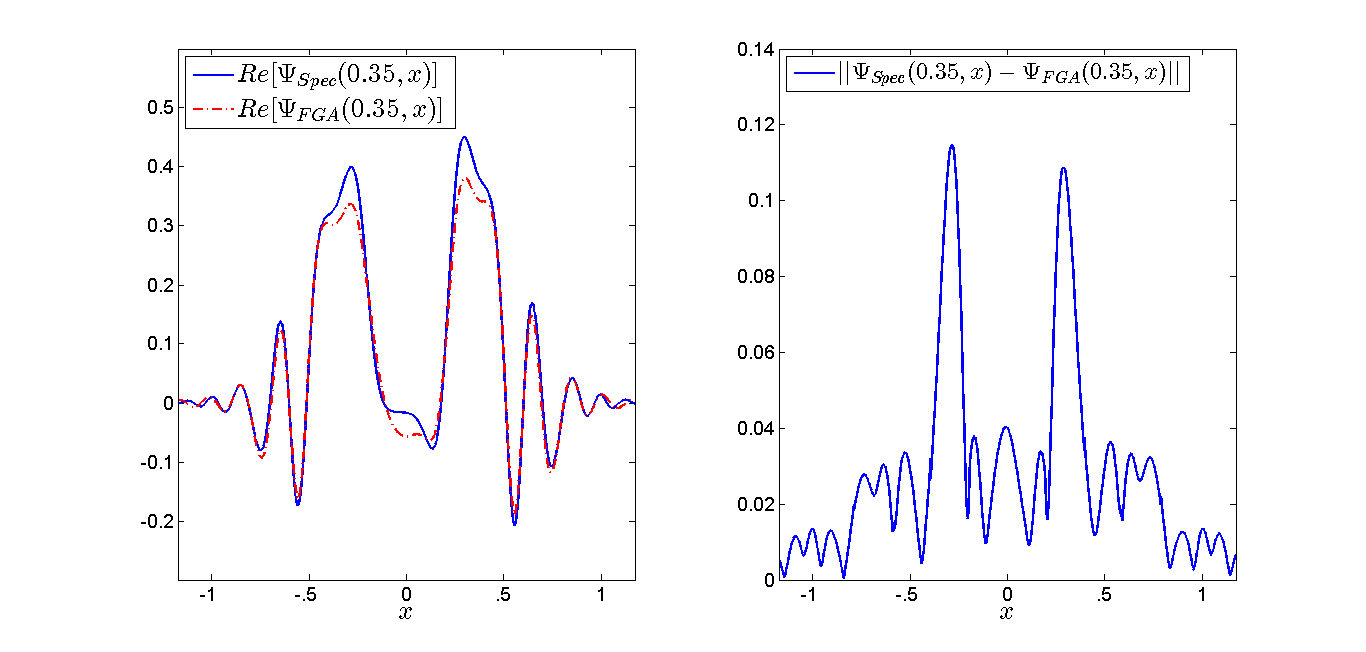}
\caption{Example \ref{ex4} plot of real parts of $\psi_{\FGA}^{\veps}(0.35,{x})$ and $\psi_{Spec}(0.35,{x})$ along side with the $L^2$ error for $\veps=1/16$.}
\label{fig17}
\end{figure}

\begin{figure}
\includegraphics[scale=.5]{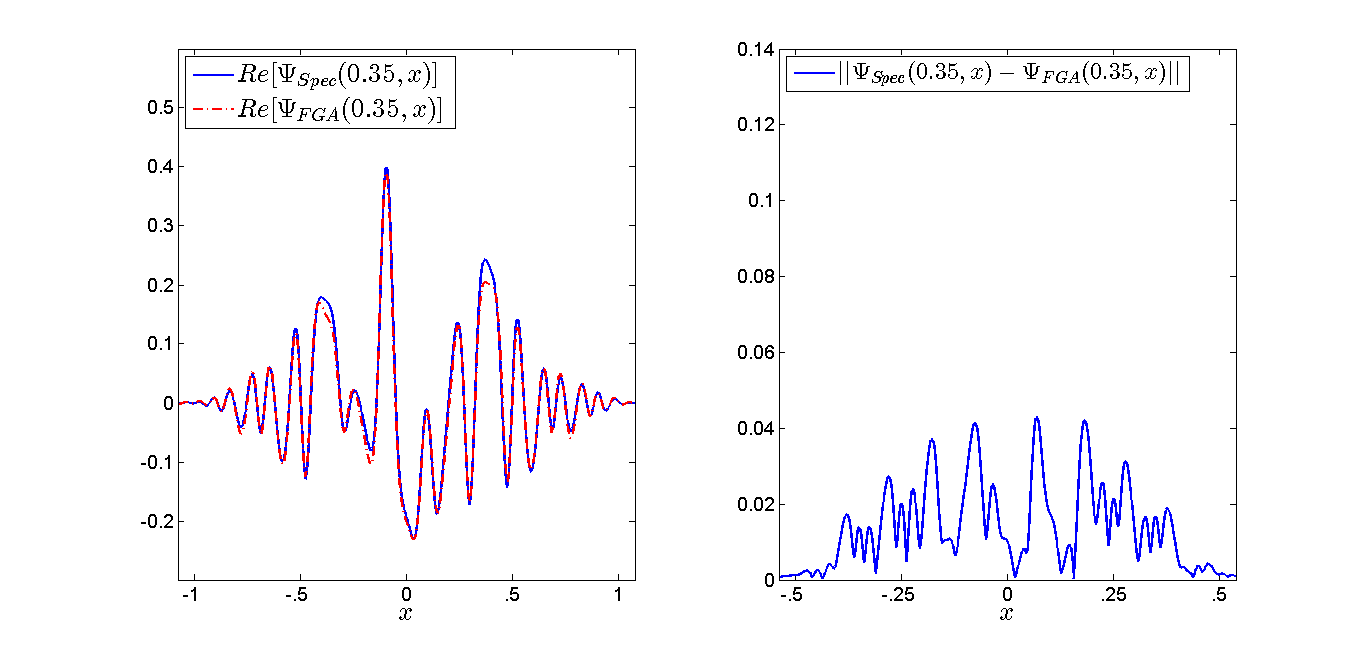}
\caption{Example \ref{ex4} plot of real parts of $\psi_{\FGA}^{\veps}(0.35,{x})$ and $\psi_{Spec}(0.35,{x})$ along side with the $L^2$ error for $\veps=1/32$.}
\label{fig18}
\end{figure}

\begin{figure}
\includegraphics[scale=.5]{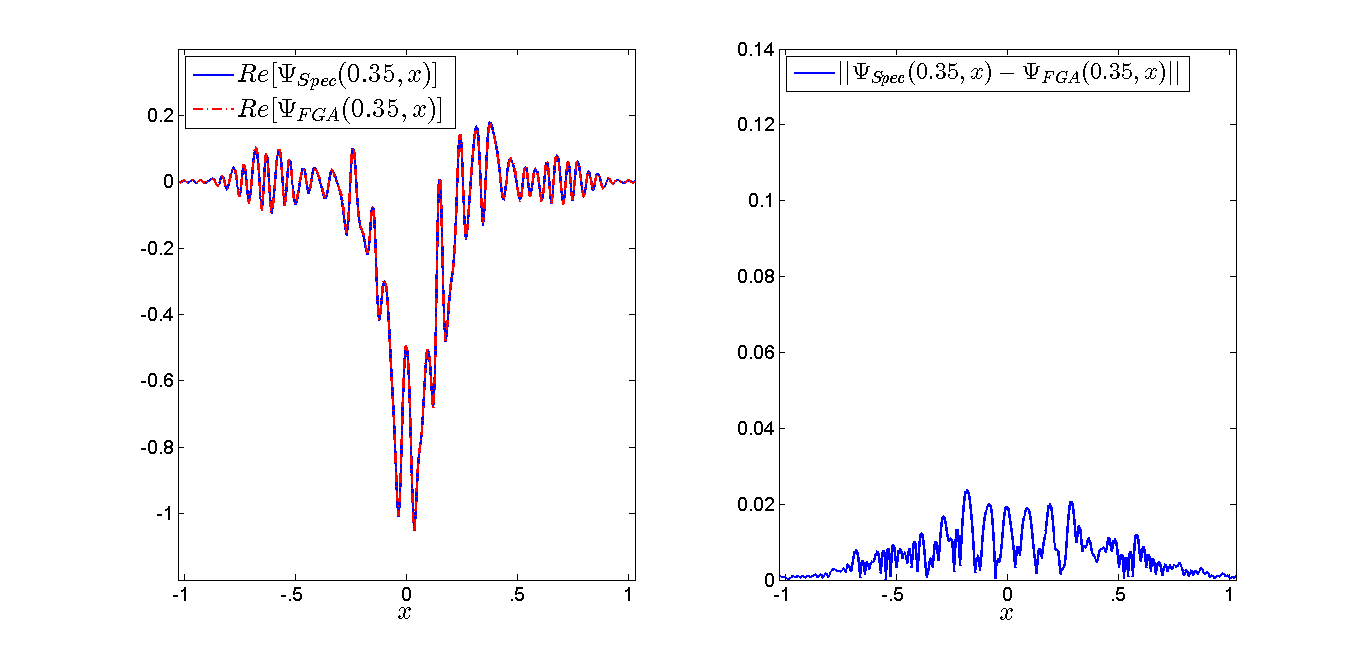}
\caption{Example \ref{ex4} plot of real parts of $\psi_{\FGA}^{\veps}(0.35,{x})$ and $\psi_{Spec}(0.35,{x})$ along side with the $L^2$ error for $\veps=1/64$.}
\label{fig19}
\end{figure}

\begin{example}\label{ex4}
In this example we choose the initial condition to be $\psi_{0}=A(x)\exp\bigl(\I S(x)/\veps\bigr)$ with $A(x)=\exp\left(-50x^{2}\right)$ and $S(x)=0.3+0.1 \sin(x-0.5)$. The exact solution is computed using the Strang Splitting spectral method. This is done at time $T=0.35$. The lattice potential used is $V_{\Gamma}(x)=\cos(x)$. We record the data in Table~\ref{ex4table}. The convergence order of the data in table \ref{ex4table} is 1.0366. We display plots of the solution for $\veps=1/8, 1/16,1/32$ and $1/64$ in Figures~\ref{fig16}, \ref{fig17}, \ref{fig18}, and \ref{fig19}.
\end{example}

\begin{table}
\begin{tabular}{|>{\columncolor[gray]{0.8}}p{4cm}|>{\columncolor[gray]{0.8}}p{4cm}|>{\columncolor[gray]{0.8}}p{4cm}|}
\hline
& Error $||\psi_{Spec}-\psi_{\FGA}^{\veps}||_{L^2}$ & Rate of Convergence\\
\hline
$\veps=1/8$  & 0.09112 &  \\
\hline
$\veps=1/16$ & 0.048907 & 0.8977\\
\hline
$\veps=1/32$ & 0.022603 & 1.1135\\
\hline
$\veps=1/64$ & 0.010555 & 1.0986\\
\hline
\end{tabular}

\smallskip
\caption{$L^2$ error of $\psi_{Spec}(0.35,{x})-\psi_{\FGA}^{\veps}(0.35,{x})$ for various values of $\veps$. The summation in $\psi_{\FGA}^{\veps}$ is over the first 8 lowest energy bands.}
\label{ex4table}
\end{table}

In the next example, we will choose initial condition projected onto one Bloch band. With this choice of initial condition, there will be no initial error.

\begin{figure}
\includegraphics[scale=.45]{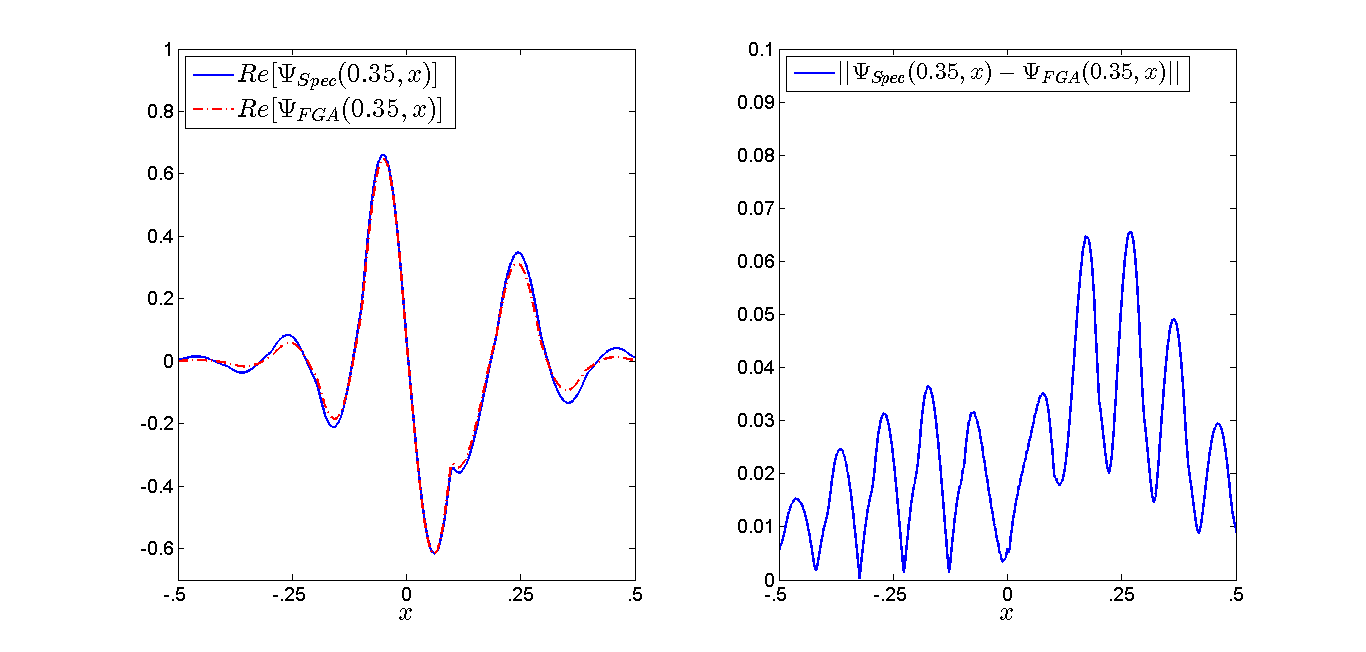}
\caption{Example \ref{ex5} plot of the real part of $\psi_{Spec}(0.35,{x})$ and $\psi_{FGA}^{\veps}(0.35,{x})$ alongside with the $L^2$ error of $\psi_{Spec}(0.35,{x})-\psi_{FGA}^{\veps}(0.35,{x})$ for example \ref{ex5}. We use $\veps=1/64$.}\label{figex5a}
\end{figure}

\begin{figure}
\includegraphics[scale=.45]{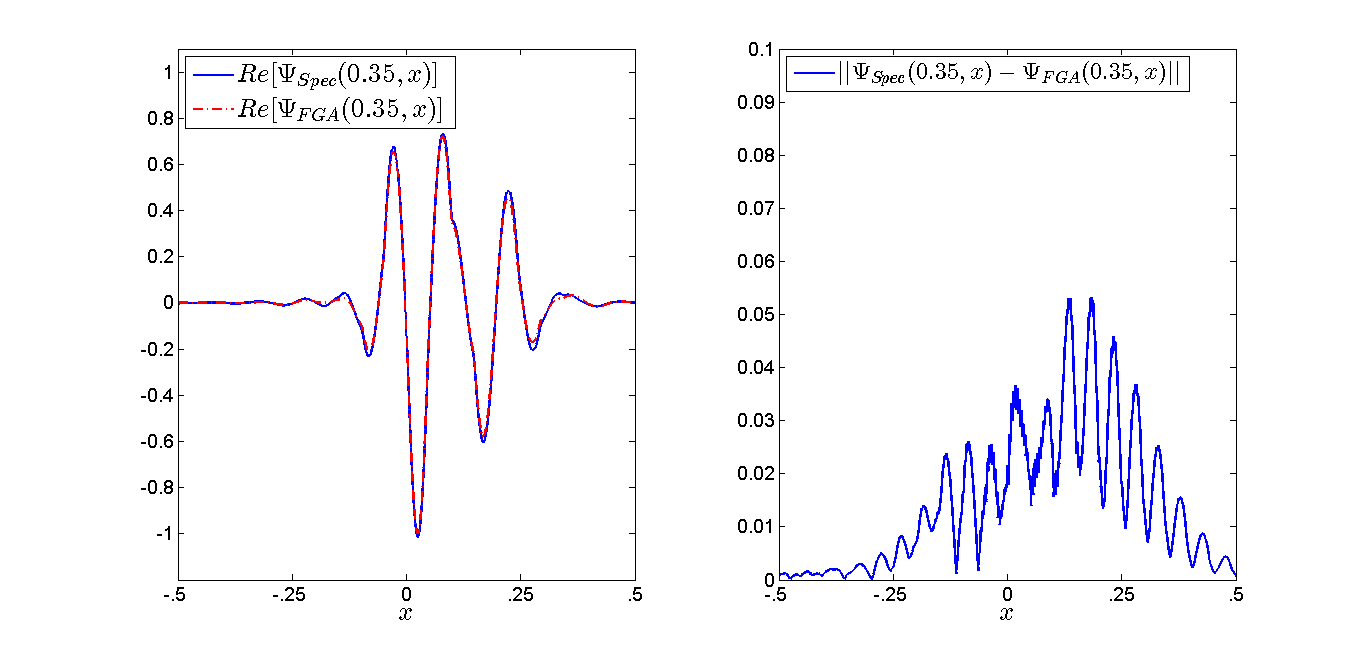}
\caption{Example \ref{ex5} plot of the real part of $\psi_{Spec}(0.35,{x})$ and $\psi_{FGA}^{\veps}(0.35,{x})$ alongside with the $L^2$ error of $\psi_{Spec}(0.35,{x})-\psi_{FGA}^{\veps}(0.35,{x})$ for example \ref{ex5}. We use $\veps=1/128$ .}\label{figex5b}
\end{figure}

\begin{figure}
\includegraphics[scale=.45]{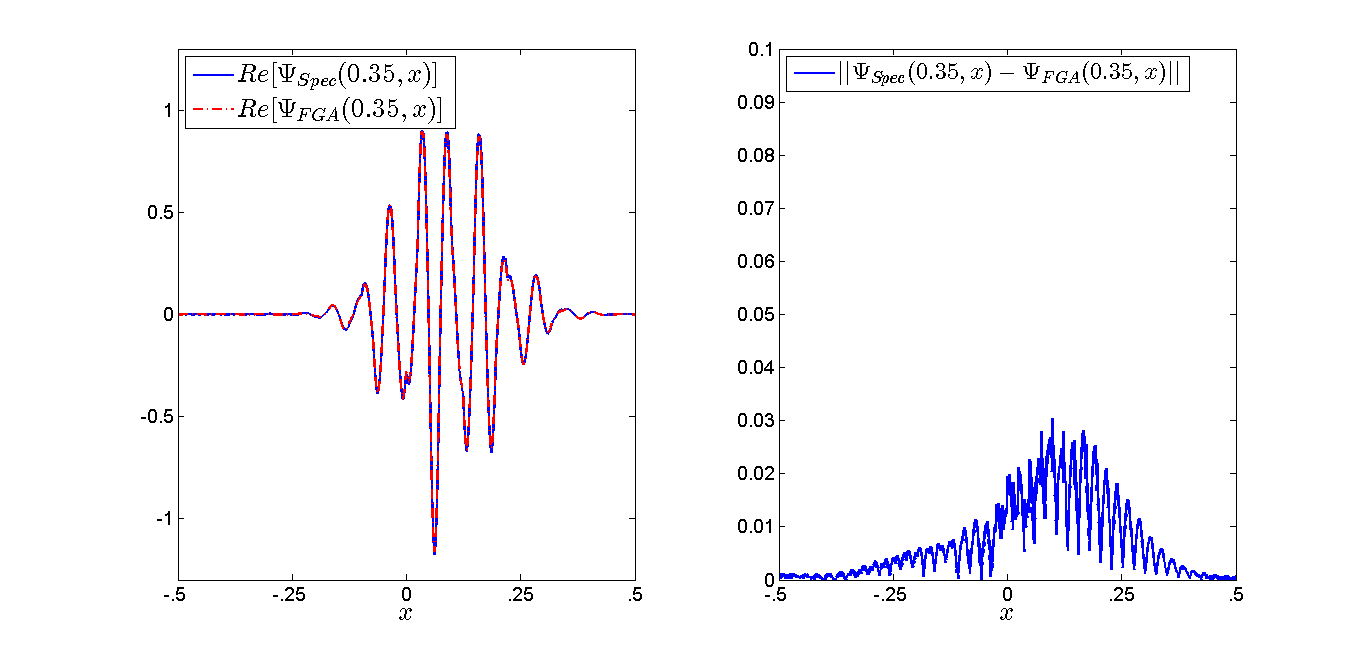}
\caption{Example \ref{ex5} plot of the real part of $\psi_{Spec}(0.35,{x})$ and $\psi_{FGA}^{\veps}(0.35,{x})$ alongside with the $L^2$ error of $\psi_{Spec}(0.35,{x})-\psi_{FGA}^{\veps}(0.35,{x})$ for example \ref{ex5}. We use $\veps=1/256$ .}\label{figex5c}
\end{figure}

\begin{example}\label{ex5}
In this example we will choose an initial condition $\Pi_{n=1}^{\mathcal{W},\veps}\psi_{0}(x)$ given by \eqref{eq:projection}  with $\psi_{0}(x)=A(x)\exp\bigl(\I S(x)/\veps\bigr)$ where $A(x)=\exp(-50x^{2})$ and $S(x)=0.3x+0.1\sin(x-0.5)$ with lattice potential $\exp(-20x^{2})$ and external potential $U(x)=0$. We compute the solution at time $T=0.35$ using the Strang Splitting spectral method and GIFGA. The $L^2$ errors are recorded in Table~\ref{ex5table}. The convergence order is 0.9814. We display plots of the solution for $\veps=1/64, 1/128$ and $1/256$ in Figures~\ref{figex5a}, \ref{figex5b}, and \ref{figex5c}.
\end{example}

\begin{table}
\begin{tabular}{|>{\columncolor[gray]{0.8}}p{4cm}|>{\columncolor[gray]{0.8}}p{4cm}|>{\columncolor[gray]{0.8}}p{4cm}|}
\hline
& Error $||\psi_{Spec}-\psi_{FGA}^{\veps}||_{L^2}$ & Rate of convergence\\
\hline
$\veps=1/64$  &  0.0269 &\\
\hline
$\veps=1/128$ &  0.0144 & 0.9015\\
\hline
$\veps=1/256$ & 0.0069 & 1.0614\\
\hline
\end{tabular}

\smallskip
\caption{$L^2$ error of $\psi_{Spec}(0.35,{x})-\psi_{FGA}^{\veps}(0.35,{x})$ for initial condition projected onto the first Bloch band.}
\label{ex5table}
\end{table}

\begin{figure}
\includegraphics[scale=.5]{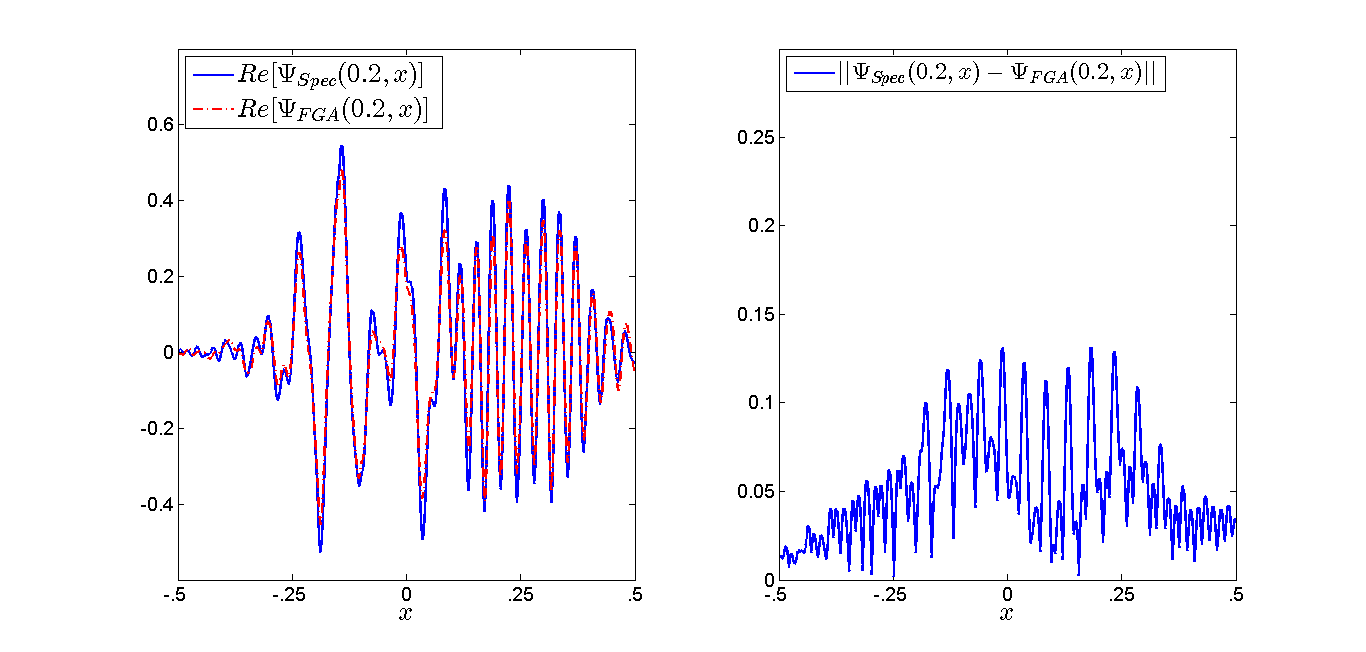}
\caption{Example \ref{ex6} plot of the real parts of $\psi_{FGA}^{\veps}(0.2,{x})$ and $\psi_{Spec}(0.2,{x})$ along side with the $L^2$ error for $\veps=1/128$.}
\label{fig10}
\end{figure}

\begin{figure}
\includegraphics[scale=.5]{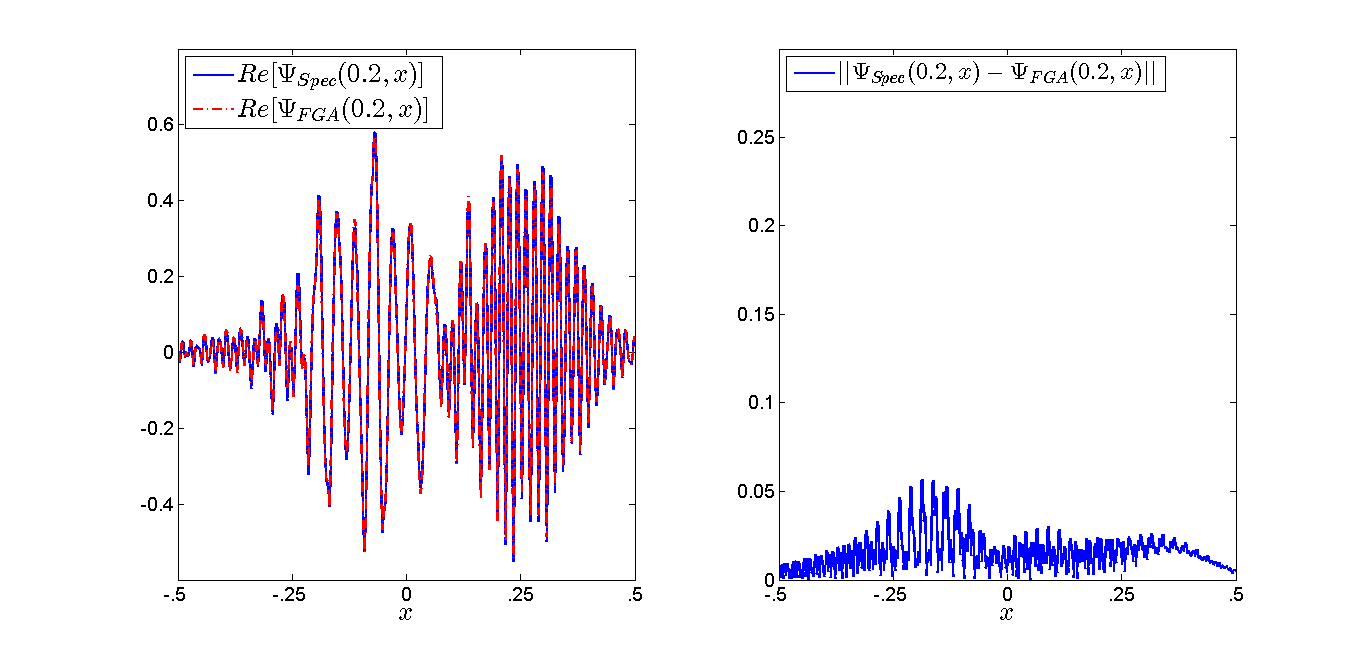}
\caption{Example \ref{ex6} plot of the real parts of $\psi_{FGA}^{\veps}(0.2,{x})$ and $\psi_{Spec}(0.2,{x})$ along side with the $L^2$ error for $\veps=1/256$.}
\label{fig11}
\end{figure}

\begin{figure}
\includegraphics[scale=.5]{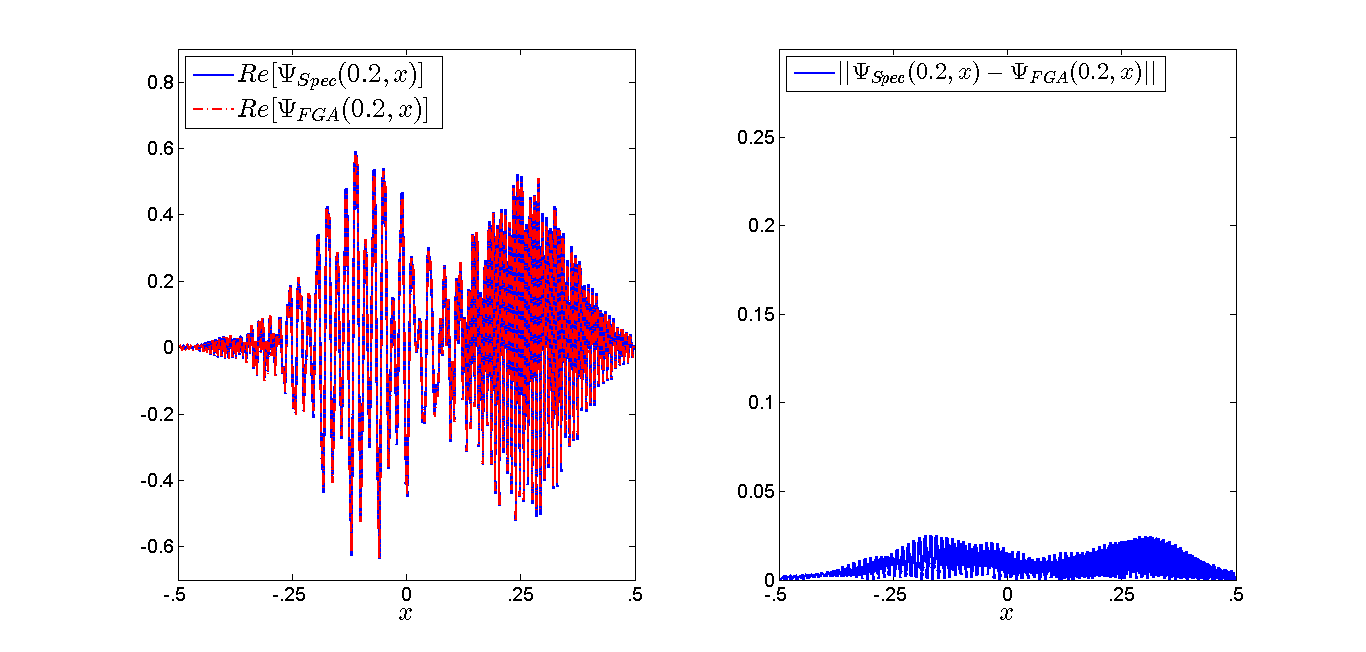}
\caption{Example \ref{ex6} plot of the real parts of $\psi_{FGA}^{\veps}(0.2,{x})$ and $\psi_{Spec}(0.2,{x})$ along side with the $L^2$ error for $\veps=1/512$.}
\label{fig12}
\end{figure}

\begin{example}\label{ex6}
In this example we choose the initial condition to be $\psi_{0}=A(x)\exp\bigl(\I S(x)/\veps\bigr)$ with $A(x)=\exp\left(-50x^{2}\right)\cos((x-0.5)/\veps)$ and $S(x)=0.3(x-0.5)+0.1\sin(x-0.5)$. The exact solution is computed using the Strang Splitting spectral method. This is done at time $T=0.2$. The potential used is $V_{\Gamma}(x)=\exp(-25x^{2})$ with external potential $U(x)=\dfrac{1}{2}x^{2}$. Our results are shown in Table~\ref{ex6table}. The convergence order of the data in table \ref{ex6table} is 0.9488. We display plots of the solution for $\veps=1/128, 1/256$ and $1/512$ in Figures~\ref{fig10}, \ref{fig11}, and \ref{fig12}.
\end{example}

\begin{table}
\begin{tabular}{|>{\columncolor[gray]{0.8}}p{4cm}|>{\columncolor[gray]{0.8}}p{4cm}|>{\columncolor[gray]{0.8}}p{4cm}|}
\hline
& Error $||\psi_{Spec}-\psi_{FGA}||_{L^2}$ & Rate of Convergence\\
\hline
$\veps=1/64$  & 0.059576 &  \\
\hline
$\veps=1/128$ & 0.038811 & .61826\\
\hline
$\veps=1/256$ & 0.015225 & 1.3500\\
\hline
$\veps=1/512$ & 0.0082833 & 0.8782\\
\hline
\end{tabular}

\smallskip
\caption{$L^2$ error of $\psi_{Spec}(0.2,{x})-\psi_{FGA}^{\veps}(0.2,{x})$ for various values of $\veps$. The summation in $\psi_{FGA}^{\veps}$ is over the first 8 lowest energy bands.}
\label{ex6table}
\end{table}


\begin{figure}
\includegraphics[scale=.5]{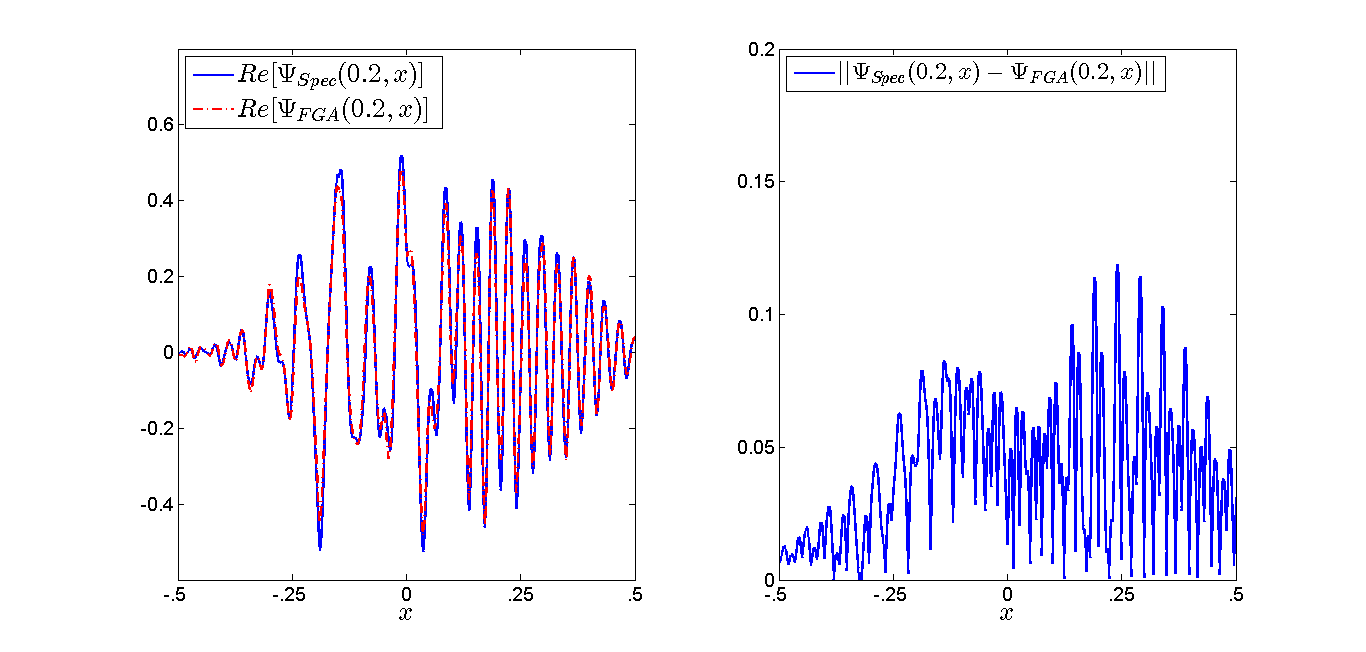}
\caption{Plot of real parts of $\psi_{FGA}^{\veps}(0.2,{x})$ and $\psi_{Spec}(0.2,{x})$ along side with the $L^2$ error for $\veps=1/128$.}
\label{fig20}
\end{figure}

\begin{figure}
\includegraphics[scale=.5]{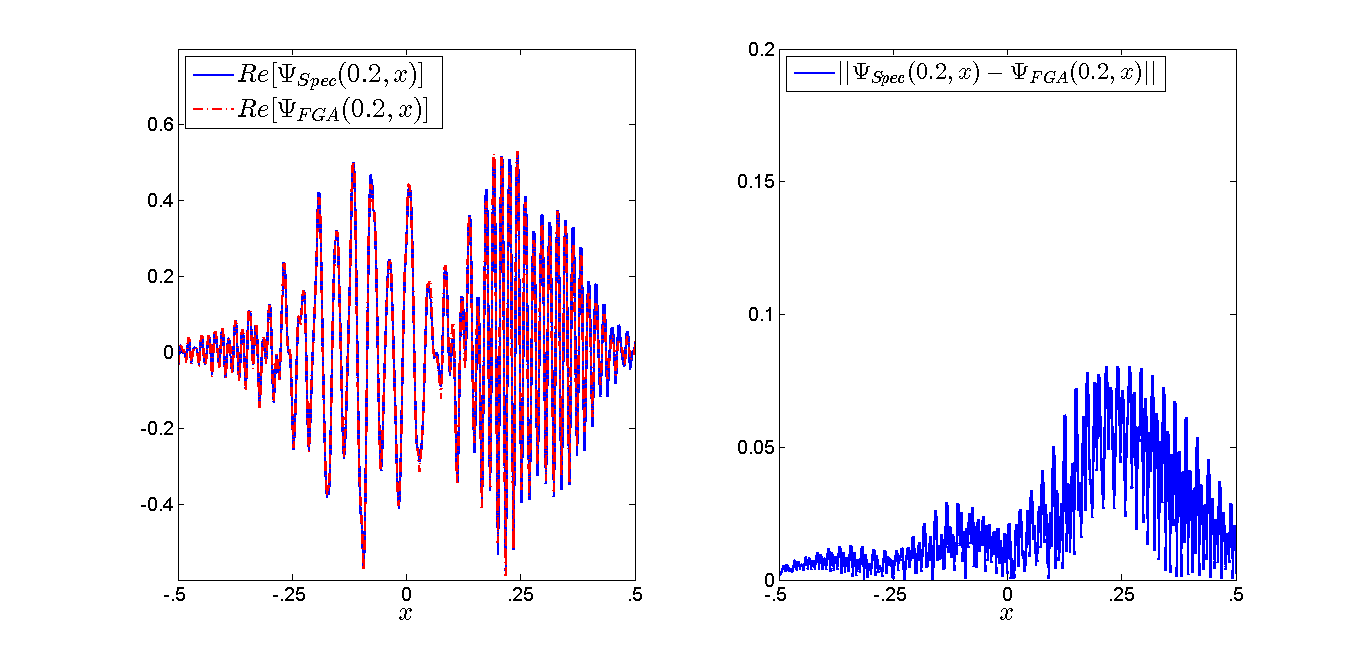}
\caption{Plot of real parts of $\psi_{FGA}^{\veps}(0.2,{x})$ and $\psi_{Spec}(0.2,{x})$ along side with the $L^2$ error for $\veps=1/256$.}
\label{fig21}
\end{figure}

\begin{figure}
\includegraphics[scale=.5]{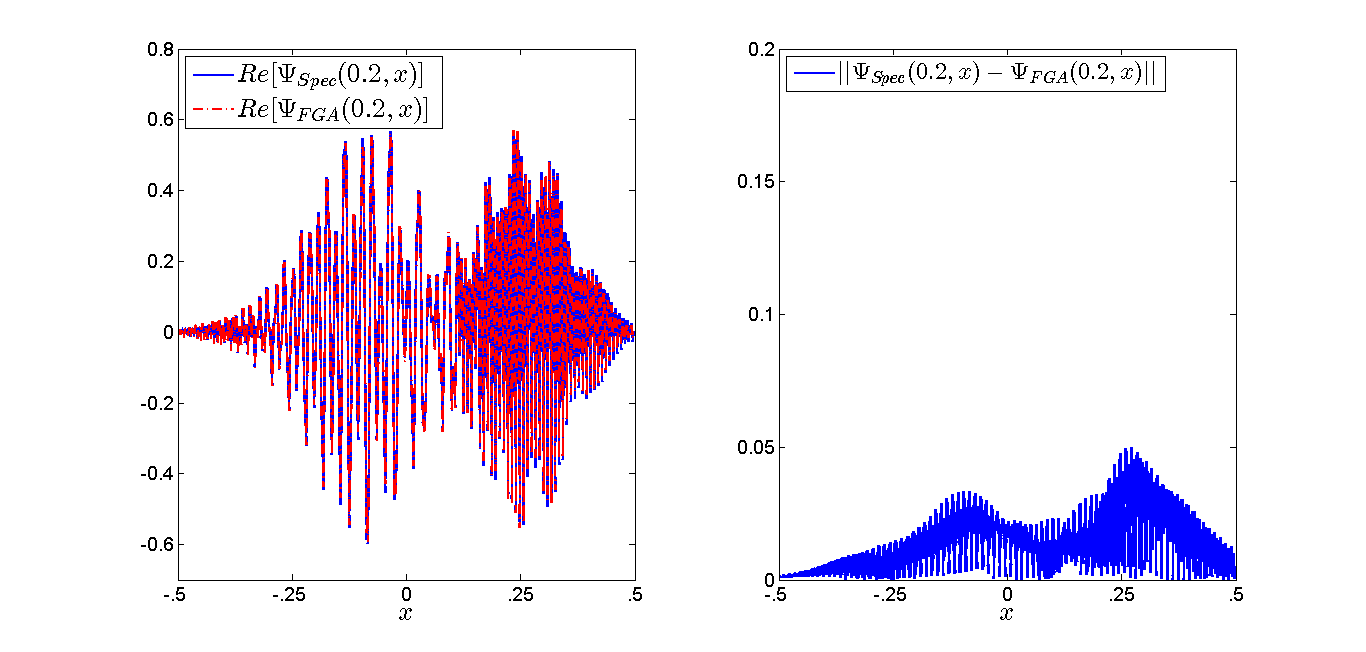}
\caption{Plot of real parts of $\psi_{FGA}^{\veps}(0.2,{x})$ and $\psi_{Spec}(0.2,{x})$ along side with the $L^2$ error for $\veps=1/512$.}
\label{fig22}
\end{figure}

\begin{example}\label{ex7}
In this example we choose the same initial condition as in Example~\ref{ex6}. All of the same parameters as in Example~\ref{ex6} will also be used. The exact solution is again computed using the Strang Splitting spectral method at time $T=0.2$. The only difference is that we change the external potential to $U(x)=\cos(x)$. The convergence order of the data in Table~\ref{ex7table} is $0.8439$. We display plots of the solution for $\veps=1/128, 1/256$ and $1/512$ in Figures~\ref{fig20}, \ref{fig21}, and \ref{fig22}.
\end{example}

\begin{table}
\begin{tabular}{|>{\columncolor[gray]{0.8}}p{4cm}|>{\columncolor[gray]{0.8}}p{4cm}|>{\columncolor[gray]{0.8}}p{4cm}|}
\hline
& Error $||\psi_{Spec}-\psi_{FGA}^{\veps}||_{L^2}$ & Rate of Convergence\\
\hline
$\veps=1/128$  & 0.039714 &  \\
\hline
$\veps=1/256$ & 0.019057 & 1.0593\\
\hline
$\veps=1/512$ & 0.012327 & 0.6285\\
\hline
\end{tabular}

\smallskip
\caption{$L^2$ error of $\psi_{Spec}(0.2,{x})-\psi_{FGA}^{\veps}(0.2,{x})$ for various values of $\veps$. The summation in $\psi_{FGA}^{\veps}$ is over the first 8 lowest energy bands.}
\label{ex7table}
\end{table}


\section{Discussion and Conclusion}\label{sec:conclusion}
  In this paper, we generalize the Herman-Kluk propagator for the linear Schr\"odinger equation (LSE),
  and develop the gauge-invariant frozen Gaussian approximation method
  for LSE with periodic potentials in the semiclassical regime. The method is invariant with respect to the gauge
  choice of the Bloch eigenfunctions, and thus avoids the numerical difficulty of computing
  gauge-dependent Berry phase. The numerical examples show that that the frozen Gaussian approximation is indeed a good approximation to the exact solution of the Schr\"odinger equation \eqref{eq:schrodinger} for $\veps\ll1$. The convergence order of our numerical results confirms the estimate given in \cite{DeLuYa:analysis}
\begin{equation}\label{eq:estimate}
||\psi_{Exact}^{\veps}(t,{x})-\psi_{FGA}^{\veps}(t,x)||_{L^{2}}  = \mathcal{O}(\veps).
\end{equation}
In future, we will study high dimensional examples where band-crossing happens quite common, and thus
requires more techniques than the scope of this paper.


\bibliographystyle{amsxport}
\bibliography{fga}

\end{document}